\documentclass[final,onefignum,onetabnum]{siamart220329}

\usepackage{braket,amsfonts}

\usepackage{array}

\usepackage[caption=false]{subfig}

\usepackage{pgfplots}

\newsiamthm{claim}{Claim}
\newsiamremark{remark}{Remark}
\newsiamremark{hypothesis}{Hypothesis}
\crefname{hypothesis}{Hypothesis}{Hypotheses}

\usepackage{algorithmic}
\usepackage{algorithm}

\usepackage{graphicx,epstopdf}
\usepackage{amsopn}

\usepackage{graphics}
\usepackage{amssymb}
\usepackage{amsmath}
\usepackage{xcolor}
\usepackage{todonotes}
\usepackage{mathrsfs}
\usepackage{multicol}
\usepackage{multirow}
\usepackage{float}
\usepackage{booktabs}
\usepackage{bm}
\usepackage{dsfont}
\usepackage{hyperref}
\hypersetup{
  breaklinks,
  colorlinks,
  citecolor={green!50!black},
  urlcolor={red!50!black},
  linkcolor={green!50!black}
}
\usepackage[nameinlink]{cleveref}
\usepackage{comment}

\makeatletter
\def\cahierline{\gdef\@cahierline}
\cahierline{\textcolor{gray}{Cahier du GERAD G-\number\year-\cahiernumber}}
\def\ps@firstpage{
\def\@oddfoot{\normalfont\scriptsize \hfil\thepage\hfil}%
  \let\@evenfoot\@oddfoot
  \def\@oddhead{\hfil\normalfont\small\@cahierline}%
  \let\@evenhead\@oddhead 
}
\def\ps@myheadings{%
    \def\@oddfoot{\normalfont\scriptsize\ttfamily\@cahierline\hfil}
    \def\@evenfoot{\normalfont\scriptsize\ttfamily\hfil\@cahierline}    \def\@evenhead{\rlap{\thepage}\hfil\upshape\footnotesize\leftmark\hfil}
    \def\@oddhead{\footnotesize\hfil{\upshape\footnotesize\rightmark}\hfil\llap{\thepage}}
    \let\@mkboth\@gobbletwo
    \let\sectionmark\@gobble
    \let\subsectionmark\@gobble
    }
\makeatother

\newcommand{\cahiernumber}{35}

\crefname{subsection}{section}{sections}
\Crefname{algorithm}{Algorithm}{Algorithms}
\newcommand{\R}{\mathds{R}}
\newcommand{\C}{\mathds{C}}

\Crefname{ALC@unique}{Line}{Lines}

\pgfplotsset{compat=1.17}

\crefname{thm}{Theorem}{Theorems}

\newtheorem{lem}{Lemma}[section]
\crefname{lem}{Lemma}{Lemmas}

\crefname{prop}{Proposition}{Propositions}

\newcommand{\pmat}[1]{\begin{pmatrix}#1\end{pmatrix}} 

\newcommand{\T}{^T\!}

\newcommand{\calT}{{\cal T}}
\newcommand{\calD}{{\cal D}}
\newcommand{\calL}{{\cal L}}
\newcommand{\calU}{{\cal U}}
\newcommand{\calA}{\cal A}
\newcommand{\calP}{\cal P}

\title{On GSOR, the Generalized Successive Overrelaxation Method for Double Saddle-Point Problems}

\author{Na Huang%
  \thanks{Department of Applied Mathematics, College of Science,
  China Agricultural University,
    Beijing, China.
    E-mail: hna@cau.edu.cn.
    Research partially supported by National Natural Science Foundation of China (No.\,12001531).
  }
  \and Yu-Hong Dai%
  \thanks{LSEC, Academy of Mathematics and Systems Science,
    Chinese Academy of Sciences, Beijing, China.
    E-mail: dyh@lsec.cc.ac.cn.
  }
  \and Dominique Orban%
  \thanks{GERAD and
    Department of Mathematics and Industrial Engineering,
    Polytechnique Montr\'eal, QC, Canada.
    E-mail: dominique.orban@gerad.ca.
    Research partially supported by an NSERC Discovery Grant.
  }
  \and Michael A. Saunders%
  \thanks{Systems Optimization Laboratory,
    Department of Management Science and Engineering,
    Stanford University, Stanford, CA, USA.
    E-mail: saunders@stanford.edu.
  }
}

\begin{document}

\maketitle

\thispagestyle{firstpage}
\pagestyle{myheadings}

\begin{abstract}
We consider the generalized successive overrelaxation (GSOR) method for solving a class of block three-by-three saddle-point problems. Based on the necessary and sufficient conditions for all roots of a real cubic polynomial to have modulus less than one, we derive convergence results under reasonable assumptions. We also analyze a class of block lower triangular preconditioners induced from GSOR and derive explicit and sharp spectral bounds for the preconditioned matrices. We report numerical experiments on test problems from the liquid crystal director model and the coupled Stokes-Darcy flow, demonstrating the usefulness of GSOR.
\end{abstract}

\begin{keywords}
   iterative methods, double saddle-point systems, saddle-point problems, matrix splitting, successive overrelaxation, preconditioning.
\end{keywords}

\begin{AMS}
   65F10, 65F50.
\end{AMS}

\section{Introduction}

We consider the double saddle-point problem
\begin{equation}\label{a1}
  {\calA}\,w:=
  \pmat{ A & B\T & C\T 
      \\ B & 0   & 0 
      \\ C & 0   & -D}
  \pmat{x \\ y \\ z} = 
  \pmat{f \\ g \\ h} =:b,
\end{equation}
where $A\in \R^{n\times n}$ and $D\in \R^{p\times p}$ are symmetric positive definite (SPD) matrices, $B\in \R^{m\times n}$ has full row rank, and $C\in \R^{p\times n}$. Linear systems like~\eqref{a1} arise from many practical applications, such as mixed and mixed-hybrid finite element approximation of the liquid crystal director model \cite{Ramage2013} and coupled Stokes-Darcy flow \cite{Cai2009,Holter20201,Holter20202,Ellingsrud2015}, and interior methods for quadratic programming problems \cite{Wright1997,ghannad-orban-saunders-2021,b32}. We emphasize that \eqref{a1} is importantly different from the block $3 \times 3$ systems considered by Huang et al.~in \cite{b33,b34}.

In principle,~\eqref{a1} can be treated as the block \(2\)\(\times\)\(2\) saddle-point problem
\begin{equation}\label{gs1}
  \pmat{ H & E\T
      \\ E &-W}
  \pmat{x \\ y}
 =\pmat{f \\ g},
  \end{equation}
which has been studied for decades \cite{Benzi2005}.
We focus here on splitting iterative methods for~\eqref{a1} by fully utilizing the special structure of $\calA$.
The generalized successive overrelaxation (GSOR) method of Bai et al.~\cite{Bai2005} is for~\eqref{gs1} with $W = 0$. We extend GSOR to~\eqref{a1} by introducing three parameters. The convergence analysis of this new GSOR method is quite different from that of stationary methods; we derive convergence conditions based on the necessary and sufficient conditions for all roots of a real cubic polynomial to have modulus less than one. We also analyze a class of block lower triangular preconditioners induced from GSOR and show that all eigenvalues of the preconditioned matrices are positive real and can be clustered by appropriate selections of parameters.

For linear systems discretized from a mixed Stokes-Darcy model, Cai et al.~\cite{Cai2009} proposed preconditioning techniques by treating~\eqref{a1} as system~\eqref{gs1} with $W=0$.
Ramage and Gartland Jr.~\cite{Ramage2013} studied a preconditioned nullspace method for solving systems~\eqref{a1} that arise from discretizations of continuum models for the orientational properties of liquid crystals, in which they also partitioned $\calA$ into a block \(2\)\(\times\)\(2\) form. Recently, based on the special structure of $\calA$ in~\eqref{a1}, several preconditioners were proposed to accelerate Krylov subspace methods. Beik and Benzi \cite{Beik20182,Beik2018} analyzed several block diagonal and block triangular preconditioners and derived bounds for the eigenvalues of the preconditioned matrices. An alternating positive semidefinite splitting (APSS) preconditioner and its relaxed variant were proposed by Liang and Zhang \cite{Liang2019} to solve double saddle-point problems arising from liquid crystal director models. The improved APSS preconditioner of Ren et al.~\cite{ren2022} and the two-parameter block triangular preconditioner of Zhu et al.~\cite{zhu2022} were also constructed to deal with the same saddle-point problem.
However, the latter preconditioners either do not fully exploit the special structure of $\calA$ or need to solve several complicated and dense linear systems at each iteration.

It is generally difficult to analyze the spectral properties of a ``full" block three-by-three matrix; i.e., one that cannot be reduced to a block \(2\)\(\times\)\(2\) matrix. Little literature exists on iterative schemes for~\eqref{a1}. Uzawa-like methods based on the splitting
\begin{equation}\label{uzawasplit}
{\calA} = 
  \pmat{A & 0 & 0 
     \\ B & -\frac{1}{\alpha}Q & 0 
     \\ C & 0 & M} -
  \pmat{0 & -B\T & -C\T 
     \\ 0 &  -\frac{1}{\alpha}Q & 0 
     \\ 0 & 0 & N}
\end{equation}
were studied by Benzi and Beik \cite{Benzi2019}, where $\alpha>0$ and the SPD matrix $Q$ are given, and $D=N-M$ with $M$ negative definite. In addition, given a parameter $\omega\neq0$, they split $\calA$ into
$$
{\calA} = \frac{1}{\omega}
  \pmat{A & B\T & 0
     \\ B & 0 & 0 
     \\ \omega C & 0 & -D}
  - \frac{1}{\omega}
  \pmat{(1-\omega)A & (1-\omega)B\T & -\omega C\T 
     \\ (1-\omega)B &  0  & 0 
     \\ 0 & 0 & -(1-\omega)D}
$$
and proposed a generalized block successive overrelaxation (GBSOR) method. The convergence analysis of these two methods is similar to that of stationary iterative schemes for block \(2\)\(\times\)\(2\) linear systems, where convergence conditions are derived from a quadratic polynomial equation of the eigenvalues of the iteration matrix. Moreover, GBSOR needs to solve four linear systems at each step: two of the form $Ax=r_1$, one $BA^{-1}B\T y=r_2$, and one $Dz=r_3$. By partitioning $\calA$ into system~\eqref{gs1} with $H=A$, Dou and Liang \cite{dou2022} construct a class of block alternating splitting implicit (BASI) iteration methods. At each step, BASI needs to solve several linear systems of the form $\alpha I + A + a B\T B + b C\T C $ and $\alpha I + D + c CC\T$, where $I$ is the identity and $\alpha$, $a$, $b$, $c$ are real scalar constants.

The paper is organized as follows. In \cref{sec:GSOR}, we present the generalized successive overrelaxation method. In \cref{sec:convergence}, convergence of GSOR is established under reasonable assumptions. Preconditioners are analyzed in \cref{sec:preconditioner}. Numerical experiments are reported in \cref{sec:numres}. Conclusions are summarized in \cref{sec:con}.

\subsection*{Notation}
For any $S\in \R^{r\times r}$, its spectral radius, inverse and transpose are denoted $\rho(S)$, $S^{-1}$ and $S^T$, respectively. For any $s\in \C^r$, its conjugate transpose is denoted $s^*$.

\section{The generalized successive overrelaxation (GSOR)  method}\label{sec:GSOR}

In this section, we present GSOR for solving the double saddle-point problem~\eqref{a1}. We consider the equivalent unsymmetric system
\begin{equation}\label{ab1}
  \widehat{{\calA}}\,w:=
  \pmat{ A & B\T & C\T 
     \\ -B & 0 & 0 
     \\ -C & 0 & D}
  \pmat{x \\ y \\ z} = 
  \pmat{f \\-g \\-h} =: \widehat{b}.
\end{equation}
Although $\widehat{\calA}$ is unsymmetric, it has certain desirable properties:
\begin{enumerate}
\item $\widehat{\calA}$ is semipositive real: $v\T \widehat{\calA} v \ge 0$ for all $v\in \R^{n+m+p}$.

\item $\widehat{\calA}$ is positive semistable; i.e., its eigenvalues have nonnegative real part.
\end{enumerate}
These properties enable convergence of the classical successive overrelaxation (SOR) method \cite{Young2014}. To improve efficiency, we modify the classical SOR method and propose a generalized version that extends the GSOR method considered in \cite{Bai2005}.

By introducing the three matrices
\begin{equation}\label{a2}
 {\calD} =
 \pmat{ A & 0 & 0 
     \\ 0 & P & 0 
     \\ 0 & 0 & D}, \quad
 {\calL} = 
 \pmat{ 0 & 0 & 0 
     \\ B & 0 & 0 
     \\ C & 0 & 0}, \quad
 {\calU} = 
 \pmat{ 0 & -B\T & -C\T 
     \\ 0 & P & 0 
     \\ 0 & 0 & 0},
\end{equation}
where $P\in\R^{m\times m}$ is SPD, 
we can split $\widehat{\calA}$ as
$$
\widehat{\calA} = {\calD}-{\calL}-{\calU}.
$$
Let $\omega$, $\tau$, and $\theta$ be three nonzero reals, $I_n$, $I_m$, and $I_p$ be identity matrices of appropriate order, and
$$
  \Omega =
  \pmat{\omega I_n & 0 & 0 
     \\ 0 & \tau I_m & 0 
     \\ 0 & 0 & \theta I_p}.
$$
Consider the following iteration for~\eqref{a1}:
\begin{eqnarray}
   w_{k+1} &=&({\calD}-\Omega {\calL})^{-1}[(I-\Omega){\calD}+\Omega{\calU}] w_{k} + ({\calD}-\Omega {\calL})^{-1}\Omega \widehat{b}\label{a3}\\
  &=:&{\calT} w_{k} + ({\calD}-\Omega {\calL})^{-1}\Omega \widehat{b}.\nonumber
\end{eqnarray}
From~\eqref{a2},
\begin{align}
{\calD}-\Omega {\calL} &=
    \pmat{ A & 0 & 0 
      \\ -\tau B & P & 0 
      \\ -\theta C & 0 & D
         }, \label{dla}
\\[4pt]
(I-\Omega){\calD}+\Omega{\calU} &=
    \pmat{(1-\omega)A & -\omega B\T & -\omega C\T 
        \\ 0 & P & 0 
        \\ 0 & 0 & (1-\theta)D}. \label{dlb}
\end{align}
Substituting~\eqref{dla} and~\eqref{dlb} into~\eqref{a3}, we obtain GSOR as stated in \Cref{gsor}.

\begin{algorithm}
\caption{The GSOR method}
\label{gsor}
\begin{algorithmic}[1]
\STATE{Choose $(x_0, y_0, z_0) \in \R^{n+m+p}$, $P\in \R^{m\times m}$ SPD, and $\omega$, $\tau$, $\theta > 0$.}
\FOR{$k = 0, 1, \dots$}
\STATE{Compute $(x_{k+1}, y_{k+1}, z_{k+1})$ according to the iteration
  \begin{equation}\label{gsorscheme} \tag{2.7}
    \left\{
    \begin{array}{l}
    x_{k+1} = x_k + \omega A^{-1}(f - Ax_k - B\T y_{k} - C\T z_k),\\[2pt]
    y_{k+1} = y_k + \tau P^{-1} (Bx_{k+1} - g),\\[2pt]
    z_{k+1} = z_k + \theta D^{-1}(Cx_{k+1} - Dz_k - h).
    \end{array}
  \right.
\end{equation}}
\ENDFOR
\end{algorithmic}
\end{algorithm}

At each step, GSOR needs to solve only three SPD systems
(of order $n$, $m$, and $p$). This is easier than in GBSOR \cite{Benzi2019}, which solves four linear systems involving $A$, $A$,  $BA^{-1}B\T$, and $D$.

Iteration scheme~\eqref{gsorscheme} can also be deduced from the splitting
\begin{equation}\label{splitting} \tag{2.6}
{\calA} = {\cal M}-{\cal N}:=
          \pmat{ \frac{1}{\omega}A & 0 & 0 
               \\ B & -\frac{1}{\tau}P & 0 
               \\ C & 0 &-\frac{1}{\theta}D
    }-
  \pmat{(\frac{1}{\omega}-1)A & -B\T &- C\T 
     \\  0 &  -\frac{1}{\tau}P & 0 
     \\  0 & 0 &(1-\frac{1}{\theta})D
       }.
\end{equation}
Therefore, GSOR is a splitting method. In particular if $\omega=1$, GSOR reduces to the Uzawa-like schemes studied in \cite{Benzi2019}, where $M=-\frac{1}{\theta}D$ and $N=\left(1-\frac{1}{\theta}\right)D$ in \eqref{uzawasplit}.

\section{Convergence analysis for GSOR} \label{sec:convergence}

First, the following two lemmas give readily verifiable necessary and sufficient conditions for all roots of a real polynomial of degree two or three, respectively, to have modulus less than one.

\begin{lem}\label{l2}{\rm\cite[Theorem 1.3]{tdperef}}
Consider the second-degree polynomial equation
\begin{equation}\label{sdpe}
  \lambda^2+a_1\lambda+a_0=0,
\end{equation}
where $a_0$ and $a_1$ are real numbers. A necessary and sufficient condition for both roots of~\eqref{sdpe} to lie in the open disk $|\lambda|<1$ is
\begin{eqnarray*}
|a_1|<1+a_0<2.
\end{eqnarray*}
\end{lem}

\begin{lem}\label{l1}{\rm\cite[Theorem 1.4]{tdperef}}
Consider the third-degree polynomial equation
\begin{equation}\label{tdpe}
  \lambda^3+a_2\lambda^2+a_1\lambda+a_0=0,
\end{equation}
where $a_0$, $a_1$, and $a_2$ are real numbers. A necessary and sufficient condition for all roots of~\eqref{tdpe} to lie in the open disk $|\lambda|<1$ is
\begin{eqnarray*}
   |a_2+a_0|<1+a_1,\quad
 & |a_2-3a_0|<3-a_1,\quad
 & a_0^2+a_1-a_0a_2<1.
\end{eqnarray*}
\end{lem}

GSOR is convergent if and only if the spectral radius of ${\calT}$ is less than $1$; i.e., $\rho({\calT})=\rho(({\calD}-\Omega {\calL})^{-1}[(I-\Omega){\calD}+\Omega{\calU}])<1$. From this point of view, we can now study the convergence properties of GSOR.

\begin{theorem}\label{mainresult}
Assume that $A$ and $D$ are SPD, and that $B$ has full row rank.
Let the maximum eigenvalues
of $A^{-1}B\T P^{-1}B$ and $A^{-1}C\T D^{-1}C$ be $\mu_{\max}$ and $\nu_{\max}$, respectively. Then GSOR is convergent if $0 < \theta < 2$ and
\begin{eqnarray*}
  0 < \omega < \dfrac{4(2-\theta)}{(2-\theta)(2+ \tau \mu_{\max})+2 \theta\nu_{\max}},
 &\qquad &
 0 < \tau < \frac{4(\omega+\theta-\omega\theta)}{\omega\theta\mu_{\max}}.
\end{eqnarray*}
\end{theorem}

\begin{proof}
Let $\lambda$ and $v=(x,y,z)$ be an eigenvalue and eigenvector of ${\calT}$. Then
$$[(I-\Omega){\calD}+\Omega{\calU}]v = \lambda ({\calD}-\Omega {\calL})v.$$
Substituting~\eqref{dla} and~\eqref{dlb} gives
\begin{align}
  &  (1-\omega)Ax - \omega B\T y - \omega C\T z = \lambda Ax, \label{eq:eiga}\\[2pt]
  &  Py = - \lambda\tau Bx + \lambda Py, \label{eq:eigb} \\[2pt]
  &  (1-\theta)Dz = -\lambda\theta Cx + \lambda Dz.\label{eq:eigc}
\end{align}
To continue the proof, we consider sufficient conditions to guarantee $|\lambda| < 1$.

If $\lambda=1-\theta$, we have $|\lambda|<1$ if and only if $0<\theta<2$. In the following, we assume that $\lambda\neq 1-\theta$ and $0<\theta<2$.

Note that we must have $\lambda\neq1$; otherwise,~\eqref{eq:eiga}--\eqref{eq:eigc} give
\begin{eqnarray}\label{a6}
   Ax+B\T y+C\T z=0,&\quad
   Bx=0, &\quad
   Dz = Cx,
\end{eqnarray}
which imply that $0=x\T Ax+x\T B\T y+x\T C\T z=x\T Ax+z\T Dz$. Given that both $A$ and $D$ are SPD, we have $x\T Ax=0$ and $z\T Dz=0$, which further means $x=0$ and $z=0$. With $B$ having full row rank, the first equality in~\eqref{a6} gives $y=0$. This contradicts the fact that $v$ is an eigenvector.

It follows from $\lambda\neq1$, $\lambda\neq1-\theta$, and~\eqref{eq:eigb}--\eqref{eq:eigc} that
\begin{eqnarray*}
   y = \frac{\lambda\tau}{\lambda-1}P^{-1}Bx,
   &\qquad& z = \frac{\lambda\theta}{\lambda+\theta-1}D^{-1}Cx.
\end{eqnarray*}
Substituting these into~\eqref{eq:eiga}, we obtain
$$
(1-\omega)Ax - \frac{\lambda\omega\tau}{\lambda-1}B\T P^{-1}Bx -  \frac{\lambda\omega\theta}{\lambda+\theta-1}C\T D^{-1}Cx = \lambda Ax.
$$
Note that $x\neq0$; otherwise, we have $Py = 0$ and $Dz =0$, which implies $v = 0$ because $P$ and $D$ are SPD. This contradicts the fact that $v$ is an eigenvector.
Therefore, premultiplying both sides 
by $x^*/(x^*Ax)$ gives
\begin{equation}\label{a10}
 1-\omega - \frac{\lambda\omega\tau}{\lambda-1}\phi(x) -  \frac{\lambda\omega\theta}{\lambda+\theta-1}\varphi(x) = \lambda ,
\end{equation}
where
\begin{eqnarray*}
  \phi(x)=\frac{x^*B\T P^{-1}Bx}{x^*Ax},
  &\qquad& 
  \varphi(x)=\frac{x^*C\T D^{-1}Cx}{x^*Ax}.
\end{eqnarray*}

First, we consider the case $x\in{\rm null}(B)$. Clearly, $\phi(x)=0$. If $x\in{\rm null}(C)$ as well, it follows from~\eqref{a10} that $\lambda=1-\omega$. To guarantee $|\lambda|<1$, we can assume that $0<\omega<2$. If $x\notin{\rm null}(C)$, with $\phi(x)=0$,~\eqref{a10} reduces to the quadratic polynomial equation
$$
\lambda^2+\left(\omega+\theta+\omega\theta\varphi(x)-2\right)\lambda+1+\omega\theta-\omega-\theta=0.
$$
By \cref{l2}, both roots $\lambda$ of this real quadratic equation satisfy $|\lambda|<1$ if and only if
$$
|\omega+\theta+\omega\theta\varphi(x)-2|<2+\omega\theta-\omega-\theta<2.
$$
If $0<\omega<2$ and $0<\theta<2$, we have $\omega+\theta\ge 2\sqrt{\omega\theta}>\omega\theta \Rightarrow \omega+\theta-\omega\theta>0$. This, along with $\varphi(x)\ge0$ yields
\begin{eqnarray}\label{result1}
  0<\omega<\frac{2(2-\theta)}{2-\theta+\theta\varphi(x)} \le 2\,.
\end{eqnarray}

Next, we consider the case $x\notin{\rm null}(B)$. Then $\phi(x)>0$ and~\eqref{a10} can be rewritten as a cubic polynomial equation $\lambda^3 + a_2\lambda^2 + a_1\lambda + a_0 = 0$, where
\begin{align*}
   a_2 &= \theta+\omega+\omega\tau\phi(x)+\omega\theta\varphi(x)-3,
\\ a_1 &= 3+\omega\theta-2\omega-2\theta-\omega\tau(1-\theta)\phi(x)-\omega\theta\varphi(x),
\\ a_0 &= \omega+\theta -\omega\theta-1.
\end{align*}
By \cref{l1}, all roots $\lambda$ of the above real cubic equation satisfy $|\lambda|<1$ if and only if
\begin{small}
\begin{align}
 \big|\, 2\omega+2\theta-\omega\theta+\omega\tau\phi(x)+\omega\theta\varphi(x)-4 \,\big|
 &< 4+\omega\theta-2\omega-2\theta-\omega\tau(1-\theta)\phi(x)-\omega\theta\varphi(x),\hspace*{-30pt}
  \label{b1}
 \\
 \big|\,\omega\tau\phi(x)+\omega\theta\varphi(x)-2\omega-2\theta +3\omega\theta \,\big|
 &< 2\omega+2\theta-\omega\theta+\omega\tau(1-\theta)\phi(x)+\omega\theta\varphi(x),\hspace*{-30pt}
\label{b2}
\\
 \theta(\omega\theta-\omega-\theta)\left(1+\varphi(x)\right)-\omega\tau(1-\theta)\phi(x) &< 0.
\label{b3}
\end{align}
\end{small}

Note that with $\phi(x)>0$, $\varphi(x)\ge0$, and $0<\theta<2$,~\eqref{b1} leads to $\tau>0$ and
\begin{equation}\label{b11}
 0 < \omega < \frac{4(2-\theta)}{4-2 \theta+ \tau(2-\theta)\phi(x)+2 \theta\varphi(x)} \le 2.
\end{equation}
It follows from~\eqref{b2} that
\begin{equation}\label{b21}
\omega\tau\theta\phi(x) < 4(\omega+\theta-\omega\theta),
\end{equation}
and
\begin{equation}\label{b22}
  \omega\tau(2-\theta)\phi(x)+2\omega\theta\varphi(x)+2\omega\theta > 0.
\end{equation}
If $\omega>0$, $\tau>0$ and $0<\theta<2$, as $\phi(x)>0$ and $\varphi(x)\ge0$, clearly~\eqref{b22} holds. This together with~\eqref{b21} implies that~\eqref{b2} holds if
\begin{eqnarray}\label{b23}
&&
0 < \tau < \frac{4(\omega+\theta-\omega\theta)}{\omega\theta\phi(x)}.
\end{eqnarray}
Inequality~\eqref{b3} holds if $0 < \theta \le 1$ because $\omega,\,\tau>0$. If $1<\theta<2$ and $0<\omega<2$, solving~\eqref{b3} leads to
$$
  \tau < \frac{\theta(\omega+\theta-\omega\theta)\left(1+\varphi(x)\right)}{\omega(\theta-1)\phi(x)}\,.
$$
Note that $\omega>0$, $\phi(x)>0$, $\omega+\theta-\omega\theta>0$ and $4(\theta-1) < \theta^2(1+\varphi(x))$, giving
$$
  \frac{4(\omega+\theta-\omega\theta)}{\omega\theta\phi(x)} < \frac{\theta(\omega+\theta-\omega\theta)\left(1+\varphi(x)\right)}{\omega(\theta-1)\phi(x)}.
$$
This implies that~\eqref{b3} holds under condition~\eqref{b23}.

To sum up, by combining~\eqref{result1},~\eqref{b11},~\eqref{b23} and the fact that $\tau(2-\theta)\phi(x)\ge0$, we know that $|\lambda|<1$ if
\begin{align}
 & 0 < \theta < 2, \label{eq:para-theta}\\[4pt]
 & 0 < \omega < \dfrac{4(2-\theta)}{4-2 \theta+ \tau(2-\theta)\phi(x)+2 \theta\varphi(x)}, \label{eq:para-omega}\\[4pt]
 & 0 < \tau < \frac{4(\omega+\theta-\omega\theta)}{\omega\theta\phi(x)}.\label{eq:para-tau}
\end{align}
For any $x\neq0$, we have $0 \le \phi(x) \le \mu_{\max}$ and $0 \le \varphi(x) \le \nu_{\max}$. Combining with~\eqref{eq:para-theta}--\eqref{eq:para-tau} completes the proof.
\end{proof}

\begin{remark}
We emphasize that parameters $\omega$, $\tau$ and $\theta$ can be chosen to satisfy the conditions derived in \cref{mainresult}. Indeed, as $0 < \omega < 2$ and $0 < \theta < 2$, we get
$$
\frac{4(\omega+\theta-\omega\theta)}{\omega\theta\mu_{\max}}
= \frac{4}{\theta\mu_{\max}} \left( 1+\frac{\theta}{\omega}-\theta \right)
>  \frac{4}{\theta\mu_{\max}} \left( 1+\frac{\theta}{2}-\theta \right)
= \frac{2(2-\theta)}{\theta\mu_{\max}}.
$$
Thus, we can first choose $\theta$ satisfying $0<\theta<2$, and then choose $\tau$ in the open interval $\left(0,\,\frac{2(2-\theta)}{\theta\mu_{\max}}\right)$.
Finally, we choose $\omega$ satisfying
\begin{eqnarray*}
  0 < \omega < \dfrac{4(2-\theta)}{(2-\theta)(2+ \tau \mu_{\max})+2 \theta\nu_{\max}}.
\end{eqnarray*}
\end{remark}

\begin{remark}\label{rem2}
If $\omega=1$,~\eqref{a10} can be simplified as
$$
\lambda^2+(\theta-2+\tau\phi(x)+\theta\varphi(x))\lambda+1-\theta-\tau\phi(x)+\tau\theta\phi(x)-\theta\varphi(x)=0.
$$
It follows from \cref{l2} that $|\lambda|<1$ holds if and only if
$$
|\theta-2+\tau\phi(x)+\theta\varphi(x)|<2-\theta-\tau\phi(x)+\tau\theta\phi(x)-\theta\varphi(x)<2.
$$
After some algebra, we see that GSOR with $\omega=1$ is convergent if
\begin{eqnarray*}
0<\theta<\frac{2}{1+\nu_{\max}},
&\qquad&
0<\tau<\frac{2(2-\theta-\theta\nu_{\max})}{(2-\theta)\mu_{\max}}.
\end{eqnarray*}
\end{remark}

\begin{remark} If $\omega=1$ and $\theta=1$, GSOR is the same as the Uzawa-like method studied in \cite[section 2.2]{Benzi2019}. In this case,~\eqref{a10} reduces to
$$
  \lambda^2+\left(\tau\phi(x)+\varphi(x)-1\right)\lambda-\varphi(x)=0.
$$
By \cref{l2}, we know that $|\lambda|<1$ holds if and only if
$$
  \big|\tau\phi(x)+\varphi(x)-1\big|<1-\varphi(x)<2.
$$
This implies that, for any $\tau$, GSOR diverges when $\nu_{\max}\ge1$. Therefore, for this special case, GSOR is convergent provided
\begin{eqnarray*}
  \nu_{\max}<1, &\qquad&
  0<\tau<\frac{2(1-\nu_{\max})}{\mu_{\max}}.
\end{eqnarray*}
This result is the same as \cite[Theorem 3]{Benzi2019}, which is the convergence theorem of the Uzawa-like method. However, we emphasize that the condition $\nu_{\max}<1$ is strong. In fact, as shown in \cref{subsec:msdm} below, the saddle-point problems from the mixed Stokes-Darcy model in porous media applications do not satisfy this condition. With \cref{rem2}, this shows that it is necessary to introduce another parameter.
\end{remark}

\section{The GSOR preconditioner}\label{sec:preconditioner}

We develop and analyze a class of block lower triangular preconditioners to accelerate Krylov methods for~\eqref{a1}.

The splitting in~\eqref{splitting} can induce a preconditioner $\cal M$ for~\eqref{a1}. The corresponding preconditioned matrix ${\cal M}^{-1}{\calA}$ has the form
$$
\pmat{ \omega I & \omega A^{-1}B\T & \omega A^{-1}C\T 
   \\ (\omega-1)\tau P^{-1}B & \omega\tau P^{-1}BA^{-1}B\T & \omega\tau P^{-1}BA^{-1}C\T 
   \\ (\omega-1)\theta D^{-1}C & \omega\theta D^{-1}CA^{-1}B\T & \theta I+\omega\theta D^{-1}CA^{-1}C\T
     }.
$$
When $\omega=1$, ${\cal M}^{-1}{\calA}$ has at least $n$ eigenvalues equal to $1$. As clustered eigenvalues 
are desirable, we consider the block lower triangular preconditioner
$$
{\calP} = 
   \pmat{ A & 0 & 0 
   \\ B & -\frac{1}{\tau}P & 0 
   \\ C & 0 &-\frac{1}{\theta}D
        }.
$$
When ${\calP}$ is used to precondition Krylov subspace methods, each step needs to solve three linear systems involving $A$, $P$, and $D$. This is more practical than the block preconditioners of \cite{Beik2018,Beik20182} and the (relaxed) APSS preconditioners of \cite{Liang2019}, which need to solve several dense linear systems involving matrices like $BA^{-1}B\T$, $D+CA^{-1}C^T$, $A+B\T B/\alpha$, and $D+ CC^T/\alpha$, where $\alpha$ is a
positive number.

To illustrate further the efficiency of our preconditioner ${\calP}$, we derive explicit and sharp bounds on the spectrum of the preconditioned matrix ${\calP}^{-1}{\calA}$. By direct calculations, we have
$$
{\calP}^{-1}{\calA} = 
  \pmat{I &  A^{-1}B\T & A^{-1}C\T 
     \\ 0 & \tau P^{-1}BA^{-1}B\T & \tau P^{-1}BA^{-1}C\T 
     \\ 0 & \theta D^{-1}CA^{-1}B\T& \theta I + \theta D^{-1}CA^{-1}C\T
       },
$$
which is similar to
$$
  \pmat{ I &   \hat{B}\T  &   \hat{C}\T 
      \\ 0 &  \tau \hat{B}\hat{B}\T & \tau \hat{B}\hat{C}\T 
      \\ 0 & \theta \hat{C}\hat{B}\T & \theta I+ \theta \hat{C}\hat{C}\T
        },
$$
where $\hat{B}=P^{-1/2}BA^{-1/2}$ and $\hat{C}=D^{-1/2}CA^{-1/2}$. Thus ${\calP}^{-1}{\calA}$ has eigenvalue 1 with multiplicity $n$, and the remaining eigenvalues are the same as those of
$$
 K = 
    \pmat{ \tau\hat{B}\hat{B}\T & \tau \hat{B}\hat{C}\T 
        \\ \theta\hat{C}\hat{B}\T & \theta I+\theta \hat{C}\hat{C}\T
         }.
$$
We can now establish the following theorem.

\begin{theorem}\label{eigenvalue_bound}
Assume that $A$ and $D$ are SPD, and $B$ has full row rank. Let the minimum and maximum eigenvalues
of $P^{-1}BA^{-1}B\T$ be $\mu_{\min}$ and $\mu_{\max}$. Let the maximum eigenvalue of $D^{-1}CA^{-1}C\T$ be $\nu_{\max}$. Then 
${\calP}^{-1}{\calA}$ has eigenvalue 1 with multiplicity at least $n$, and the remaining eigenvalues lie in the interval
$$
\left[\, \frac{\underline{\Lambda}-\sqrt{\underline{\Lambda}^2-4\tau\theta\mu_{\min}}}{2},\,\,
\frac{\overline{\Lambda}+\sqrt{\overline{\Lambda}^2-4\tau\theta\mu_{\max}}}{2}\,
\right],
$$
where
\begin{eqnarray}\label{def_phi}
\underline{\Lambda}=\theta(1+\nu_{\max})+\tau\mu_{\min},\quad && \quad
\overline{\Lambda} =\theta(1+\nu_{\max})+\tau\mu_{\max}.
\end{eqnarray}
\end{theorem}

\begin{proof}
We need to estimate spectral bounds for $K$. Let $\lambda$ be an eigenvalue of $K$ and $(y\T,z\T)\T$ be a corresponding eigenvector. With $\tau>0$ and $\theta>0$ we see that $K$ is similar to a symmetric matrix, and hence $\lambda$ is real.
Also,
\begin{align}
   \tau   \hat{B}\hat{B}\T y + \tau   \hat{B}\hat{C}\T z
      &= \lambda y, \label{eq:pre-eig-a}
\\ \theta \hat{C}\hat{B}\T y + \theta z
                             + \theta \hat{C}\hat{C}\T z
      &= \lambda z. \label{eq:pre-eig-b}
\end{align}
We obtain estimates of $\lambda$ by considering two cases separately.

Case I: $z\in \mathrm{null}(\hat{C}^T)$. Clearly, $\tau \hat{B}\hat{B}\T y = \lambda y$ and $ \theta \hat{C}\hat{B}\T y = (\lambda-\theta) z$. This implies that $\lambda=\theta$ or $\lambda$ is an eigenvalue of $\hat{B}\hat{B}\T$. Note that $\hat{B}\hat{B}\T=P^{-1/2}BA^{-1}B\T P^{-1/2}$ is similar to $P^{-1}BA^{-1}B\T$, so that $\lambda=\theta$ or $\tau\mu_{\min}\le\lambda\le\tau\mu_{\max}$.

Case II: $z\notin \mathrm{null}(\hat{C}^T)$. We only consider the case $\lambda\notin[\tau\mu_{\min},\,\tau\mu_{\max}]$, so that $\lambda I - \tau \hat{B}\hat{B}\T$ is nonsingular. With~\eqref{eq:pre-eig-a}, this leads to
$
  y=\tau (\lambda I - \tau \hat{B}\hat{B}^T)^{-1} \hat{B}\hat{C}\T z .
$
Substituting into~\eqref{eq:pre-eig-b} gives
\begin{eqnarray}\label{a8}
  \tau \theta \hat{C}\hat{B}\T(\lambda I - \tau \hat{B}\hat{B}^T)^{-1} \hat{B}\hat{C}\T z+\theta z+\theta \hat{C}\hat{C}\T z = \lambda z.
\end{eqnarray}
As $
\left(I-\frac{\tau}{\lambda}\hat{B}\T\hat{B}\right)^{-1}
= I + \tau \hat{B}\T(\lambda I - \tau \hat{B}\hat{B}\T)^{-1} \hat{B}
$,
\eqref{a8} yields
\begin{eqnarray}\label{ab8}
    \theta \hat{C} \left(I-\frac{\tau}{\lambda}\hat{B}\T\hat{B}\right)^{-1}\hat{C}\T z = (\lambda-\theta) z.
\end{eqnarray}
We assert that $\lambda>0$. Otherwise, we have 
$$(\lambda-\theta) z\T z <0 
  \qquad \mbox{and} \qquad
  z\T\hat{C} \left(I-\dfrac{\tau}{\lambda}\hat{B}\T\hat{B}\right)^{-1}\hat{C}\T z\ge 0,$$
which contradicts~\eqref{ab8}.

If $\lambda>\tau\mu_{\max}$, as the matrices $\hat{B}\T\hat{B}$ and $\hat{B}\hat{B}\T$ have the same nonzero eigenvalues, it holds for any $0\neq u\in \R^{n}$ that
$$
u\T\left(I-\frac{\tau}{\lambda}\hat{B}\T\hat{B}\right) u \ge  \left(1-\frac{\tau\mu_{\max}}{\lambda}  \right)u\T u > 0.
$$
With~\eqref{ab8} and the fact that $\hat{C}\hat{C}\T = D^{-1/2}CA^{-1}C\T D^{-1/2}$ is similar to $D^{-1}CA^{-1}C\T$, this leads to
$$
  \lambda-\theta \le \theta\left(1-\frac{\tau\mu_{\max}}{\lambda}  \right)^{-1} \frac{z\T\hat{C} \hat{C}\T z}{z\T z}
  \le \theta\left(1-\frac{\tau\mu_{\max}}{\lambda}  \right)^{-1}\nu_{\max}.
$$
Solving this inequality for $\lambda$ gives
$$
  \frac{\overline{\Lambda}-\sqrt{ \overline{\Lambda}^2-4\tau\theta\mu_{\max}}}{2}
  \le \lambda \le
  \frac{\overline{\Lambda}+\sqrt{\overline{\Lambda}^2-4\tau\theta\mu_{\max}}}{2}.
$$
We can directly check that
$$
\overline{\Lambda}-\sqrt{\overline{\Lambda}^2-4\tau\theta\mu_{\max}}\,
\le 2\max\{\theta, \tau\mu_{\max}\} \le
\overline{\Lambda}+\sqrt{\overline{\Lambda}^2-4\tau\theta\mu_{\max}}.
$$
Therefore, $\lambda$ admits the upper bound
\begin{equation}\label{uperbound}
\lambda \le \frac{\overline{\Lambda}+\sqrt{\overline{\Lambda}^2-4\tau\theta\mu_{\max}}}{2}.
\end{equation}

If $\lambda<\tau\mu_{\min}$, it can be verified that
$$
z\T\hat{C} \left(I-\frac{\tau}{\lambda}\hat{B}\T\hat{B}\right)^{-1}\hat{C}\T z
\ge
\left(1-\frac{\tau\mu_{\min}}{\lambda}\right)^{-1}z\T\hat{C}\hat{C}\T z
\ge \left(1-\frac{\tau\mu_{\min}}{\lambda}\right)^{-1} \nu_{\max}z\T z.
$$
Combining with~\eqref{ab8} gives
$$
\lambda-\theta \ge
\left(1-\frac{\tau\mu_{\min}}{\lambda}\right)^{-1} \theta \nu_{\max}
= \frac{\lambda\theta\nu_{\max}}{\lambda- \tau\mu_{\min}}.
$$
Note that $\lambda<\tau\mu_{\min}$ and the inequality can be simplified as
$$
\lambda^2-(\theta+\theta \nu_{\max}+\tau\mu_{\min}) \lambda+\tau\theta\mu_{\min}  \le 0.
$$
By the definition of $\underline{\Lambda}$, we have
$$
\frac{\underline{\Lambda}-\sqrt{\underline{\Lambda}^2-4\tau\theta\mu_{\min}}}{2}
\le \lambda \le
\frac{\underline{\Lambda}+\sqrt{\underline{\Lambda}^2-4\tau\theta\mu_{\min}}}{2}.
$$
Similarly, we can check that
$$
 \underline{\Lambda}-\sqrt{\underline{\Lambda}^2-4\tau\theta\mu_{\min}} \,
\le 2\min\{\theta,\tau\mu_{\min}\} \le
 \underline{\Lambda}+\sqrt{\underline{\Lambda}^2-4\tau\theta\mu_{\min}} .
$$
This implies that $\lambda$ admits the lower bound
$$
\tau\mu_{\min}>\lambda\ge\frac{\underline{\Lambda}-\sqrt{\underline{\Lambda}^2-4\tau\theta\mu_{\min}}}{2}.
$$
Combining with~\eqref{uperbound} and the bounds derived in Case I completes the proof.
\end{proof}

\begin{remark}
To precondition the equivalent unsymmetric system~\eqref{ab1}, we can use the block lower triangular preconditioner
$$
\widehat{\calP} =
   \pmat{ A & 0 & 0 
   \\ -B & \frac{1}{\tau}P & 0 
   \\ -C & 0 &\frac{1}{\theta}D
        }.
$$
Because $\widehat{\calP}^{-1}\widehat{{\calA}}= {\calP}^{-1}{\calA}$, the preconditioned matrix $\widehat{\calP}^{-1}\widehat{\calA}$ possesses the same spectral bounds as in \cref{eigenvalue_bound}.
\end{remark}

\begin{remark}
\cref{eigenvalue_bound} shows that the preconditioned matrices ${\calP}^{-1}{\calA} = \widehat{\calP}^{-1}\widehat{\calA}$ are positive stable. Moreover, their condition number is bounded by
$$
\max\left\{\,\frac{\overline{\Lambda}+\sqrt{\overline{\Lambda}^2-4\tau\theta\mu_{\max}}}{2},\, \,
\frac{\overline{\Lambda}+\sqrt{\overline{\Lambda}^2-4\tau\theta\mu_{\max}}}
{\underline{\Lambda}-\sqrt{\underline{\Lambda}^2-4\tau\theta\mu_{\min}}}
\,\right\}.
$$
Using~\eqref{def_phi}, we obtain
$$
\frac{\overline{\Lambda}+\sqrt{\overline{\Lambda}^2-4\tau\theta\mu_{\max}}}{2}
\le \overline{\Lambda}
= \theta(1+\nu_{\max})+\tau\mu_{\min}
$$
and
\begin{eqnarray*}
&&\frac{\overline{\Lambda}+\sqrt{\overline{\Lambda}^2-4\tau\theta\mu_{\max}}}
{\underline{\Lambda}-\sqrt{\underline{\Lambda}^2-4\tau\theta\mu_{\min}}}
= \frac{\left(\overline{\Lambda}+\sqrt{\overline{\Lambda}^2-4\tau\theta\mu_{\max}}\right)
\left( \underline{\Lambda}+\sqrt{\underline{\Lambda}^2-4\tau\theta\mu_{\min}}\right)}
{4\tau\theta\mu_{\min}}
\le \frac{\underline{\Lambda}\overline{\Lambda}}{\tau\theta\mu_{\min}}\\[8pt]
&& = \frac{\theta^2(1+\nu_{\max})^2+\tau^2\mu_{\min}\mu_{\max}+\tau\theta(1+\nu_{\max})(\mu_{\max}+\mu_{\min})
 }{\tau\theta\mu_{\min}}
 \\[8pt]
&& =
\frac{\theta}{\tau}\frac{(1+\nu_{\max})^2}{\mu_{\min}}+\frac{\tau}{\theta}\mu_{\max}
+(1+\nu_{\max})\left(1+\frac{\mu_{\max}}{\mu_{\min}}\right).
\end{eqnarray*}
This shows that the matrices ${\calP}^{-1}{\calA}$ and $\widehat{\calP}^{-1}\widehat{\calA}$ will be well-conditioned given appropriate selections of parameters $\tau$, $\theta$ and matrix $P$ when $\nu_{\max}$ is not too large.\footnote{This is a reasonable request. As shown in \cref{sec:numres} below, $\nu_{\max}$ of the saddle-point systems from the liquid crystal directors model and the mixed Stokes-Darcy model in porous media applications is $0.1750$ and $1.0057$, respectively.}
\end{remark}

\section{Numerical experiments}
\label{sec:numres}

We present the results of numerical tests to examine the feasibility and effectiveness of GSOR. All experiments were run using MATLAB R2015b on a PC with an Intel(R) Core(TM) i7-8550U CPU @~1.8GHz and 16GB of RAM.
The initial guess is taken to be the zero vector, and the algorithms are terminated when the number of iterations exceeds $10^5$ or
$$
 {\rm Res}:= \|b-{\calA} w_k\|_2 / \|b\| \le 10^{-8},
$$
where $w_k$ is the current approximate solution. We report the number of iterations, the CPU time, and the final value of the relative residual, denoted by “Iter”, “CPU” and “Res”, respectively.

For our GSOR method, we tried just a few values of the parameters $\omega$, $\tau$ and $\theta$.
We compared our method with the Uzawa-like method (denoted ``Uzawa") and the generalization of the block SOR method  (denoted ``GBSOR") studied in \cite[Section 2.2 and Section 3]{Benzi2019}, respectively. We emphasize that the Uzawa method is a special case of our GSOR method with $P=Q$, $\omega=1$, and $\theta=1$. For GBSOR, based on \cite[Theorem 5]{Benzi2019}, we chose $\omega = s/4,\,s/2,\,3s/4$ (denoted ``GBSORa", ``GBSORb", ``GBSORc", respectively), where $s = 2/(1+\sqrt{\nu_{\max}})$ is the upper bound of the convergence interval for the parameter $\omega$. We used the function ``eigs" to compute $\nu_{\max}$.

We also tested Krylov methods for~\eqref{a1} or~\eqref{ab1}, such as MINRES, GMRES, and BICGSTAB. For preconditioned MINRES (denoted ``BPMINRES"), we use the block diagonal preconditioner
$$
 \pmat{ A & 0  & 0 
     \\ 0 & BA^{-1}B\T & 0
     \\ 0 & 0 & D+CA^{-1}C\T
      }.
$$
For $D=0$, this block diagonal preconditioner has been studied in \cite{Beik2018}. For preconditioned GMRES, we test the GSOR preconditioner $\calP$ with $\tau=\theta=1$ (denoted ``GPGMRES") and the block triangular preconditioner \cite{Beik2018} (denoted ``BPGMRES")
$$
 \pmat{A &  B\T & C\T 
    \\ 0 & -BA^{-1}B\T  & 0
    \\ 0 &  0  & -(D+CA^{-1}C\T)
       }.
$$


\subsection{Saddle-point systems from the liquid crystal directors model}\label{subsec:lcdm}

Continuum models for the orientational properties of liquid crystals require minimization of free energy functionals of the form
\begin{equation}\label{lcdm}
  {\cal F}[u,v,w,U]=\frac{1}{2}\int_0^1[(u_z^2+v_z^2+w_z^2)-\eta^2(\beta+w^2)U_z^2]{\rm d}z,
\end{equation}
where $u$, $v$, $w$, and $U$ are functions of $z\in[0,1]$ subject to suitable end-point conditions, $u_z$, $v_z$, $w_z$, and $U_z$ denote the first derivatives of the corresponding functions with respect to $z$, and $\eta$ and $\beta$ are prescribed positive parameters. By discretizing with a uniform piecewise-linear finite element scheme with $N+1$ cells using nodal quadrature and the prescribed boundary conditions, we minimize the free energy~\eqref{lcdm} under the unit vector constraint. We apply the Lagrange multiplier method to solve this discretized minimization model, and Newton's method to solve the nonlinear equations from the first-order conditions of the Lagrangian. Each step involves the solution of a linear system of the form~\eqref{a1} with $n=3N$ and $m=p=N$. For more details, we refer to \cite{Ramage2013}.


In our numerical experiments we set $\eta=\sqrt{3}\pi/4$ and $\beta=0.5$, which is known as the critical switching value. The discretized matrix $A$ is tridiagonal, so in all algorithms we solve systems $Ax=r$ directly by the function ``\verb|\|", which uses a tridiagonal solver. We set $P=BA^{-1}B\T$ and solve systems $P y=r$ using Cholesky factorization.
Numerical results 
are listed in \Cref{tab11,tab12} with $N=1023,2047,4095,8191,16383$, where the parameter choices for GSOR and the corresponding notation are as follows:  

\begin{center}\scriptsize
   \medskip
   \begin{tabular}{|c |c |c |c |c|}
    \hline
    Method & GSORa & GSORb & GSORc & GSORd
    \\ \hline
    $(\omega,\tau,\theta)$  & $(1,1,1)$ & $(0.95,0.95,0.95)$& $(0.9,0.8,1)$ & $(0.95,1,0.95)$
    \\ \hline
   \end{tabular}
   \medskip
\end{center}
For this problem, $A-C\T D^{-1}C$ is SPD, which guarantees convergence of the Uzawa-like method \cite[Theorem 3]{Benzi2019}. We set $Q = BA^{-1}B\T$ and $\alpha = 1-\nu_{\max}$, where $\nu_{\max}=0.1750$ is the maximum eigenvalue of $A^{-1}C\T D^{-1}C$. MINRES, GMRES and BICGSTAB without preconditioning failed to solve this problem. (For GMRES, we set the restart frequency to $100$.) BICGSTAB hit an error condition.
Therefore in \cref{tab12} we only report results from preconditioned MINRES and preconditioned GMRES.

\begin{table}[htb]\scriptsize   
  \centering
  \caption{CPU time for Cholesky factorization of $P = B A^{-1} B^T$.}\label{tab11}
  \begin{tabular}{|c |c |c |c |c |c |c |c |c |}
  \hline
  $N$& 1023&2047&4095&8191&16383\\\hline
  $n$ & 3069 & 6141 & 12285 & 24573 & 49149\\\hline
  $m$ & 1023&2047&4095&8191&16383\\\hline
  $p$ & 1023&2047&4095&8191&16383\\\hline
  $n+m+p$ & 5115 &10235 &20475&40955&81915\\\hline
  $BA^{-1}B\T$ & 0.086 & 0.43 & 1.95 & 8.93& 145.2
  \\ \hline
  \end{tabular}
\end{table}

\begin{table}[p]\scriptsize   
  \centering
  \caption{Numerical results for saddle-point systems from the liquid crystal directors model.}\label{tab12}
  \begin{tabular}{|c|c|c|c|c|c|c|c|c|c|c|}
    \hline
      & $N$  & 1023     & 2047 & 4095 & 8191 &16383\\\hline
      &Iter  & 24       &   24 &  25 & 25 & 26\\
GSORa &CPU   & 0.20 & 1.08 & 4.81 & 23.86& 226.04 \\
      &Res   & 5.95e-09 & 8.42e-09 &5.44e-09 & 7.70e-09 & 4.96e-09\\\hline

      &Iter  & 15       &   15 &  16 & 16 & 16 \\
GSORb &CPU   & 0.12 & 0.64  & 3.03& 15.35 & 152.79  \\
      &Res   & 5.02e-09 & 7.09e-09 & 2.53e-09 & 3.57e-09 & 5.05e-09\\\hline

      &Iter  & 16       &  17 &  17 & 17 & 17\\
GSORc &CPU   & 0.14&  0.71& 3.29 & 15.28 & 160.58\\
      &Res   & 8.41e-09 & 1.82e-09& 2.01e-09 & 2.34e-09 & 2.88e-09\\\hline

      &Iter  & 14       &   14 &  14 & 14 & 14\\
GSORd &CPU   &0.11 &  0.59& 2.83 & 12.64 & 132.68\\
      &Res   & 7.30e-10 &9.62e-10 & 1.31e-09 & 1.81e-09 & 2.53e-09\\\hline

      &Iter  & 18       &   18 &  18 & 20  & 20\\
UZAWA &CPU   & 0.14 & 0.76& 3.51 & 18.35 & 187.64\\
      &Res   & 4.01e-09 & 5.67e-09 & 8.02e-09 & 1.56e-09  & 2.22e-09\\\hline

      &Iter  & 72       &   73 &  75 & 76 & 77\\
GBSORa &CPU   & 0.56& 3.12 & 14.80 & 70.90 & 773.60\\
      &Res   & 8.38e-09 & 9.23e-09& 7.90e-09 & 8.69e-09 & 9.57e-09\\\hline

      &Iter  & 29     &  30 &  30 & 31 & 31\\
GBSORb &CPU   & 0.23 & 1.31 & 5.98 & 28.03 &294.17\\
      &Res   & 7.91e-09 & 5.95e-09& 8.42e-09 & 6.32e-09 & 8.95e-09\\\hline

      &Iter  & 35   &   36  &  36 & 37 & 38 \\
GBSORc &CPU   & 0.27& 1.59  & 7.13 & 34.63 & 356.66\\
      &Res   & 7.59e-09& 6.34e-09& 8.97e-09 & 7.51e-09 & 6.30e-09\\\hline

      &Iter  & 13   &   13 &  13  & 14 & 14\\
BPMINRES &CPU   & 3.50 &  18.10  & 96.25 & 837.49 & 10736.63\\
      &Res   & 2.67e-09 & 4.75e-09 & 9.35e-09 & 1.41e-09 & 1.51e-09\\\hline

      &Iter  & 8  &   8  &  8 & 8 & 8\\
BPGMRES &CPU   & 3.35&  18.63  & 91.74   & 543.05 & 10609.10\\
      &Res   & 6.42e-09  & 1.34e-08 & 9.35e-09  & 7.74e-08 &1.86e-07\\\hline

      &Iter  & 8  &   8  &  8 & 8 & 8\\
GPGMRES &CPU   & 1.80&  8.17  & 40.53  & 251.52 & 3548.65 \\
      &Res   & 9.00e-10 & 1.92e-09 & 4.99e-09  & 1.66e-08 &7.79e-08 \\\hline
  \end{tabular}
\end{table}

\begin{figure}[p]   
    \centering
    \includegraphics[width=0.4\textwidth]{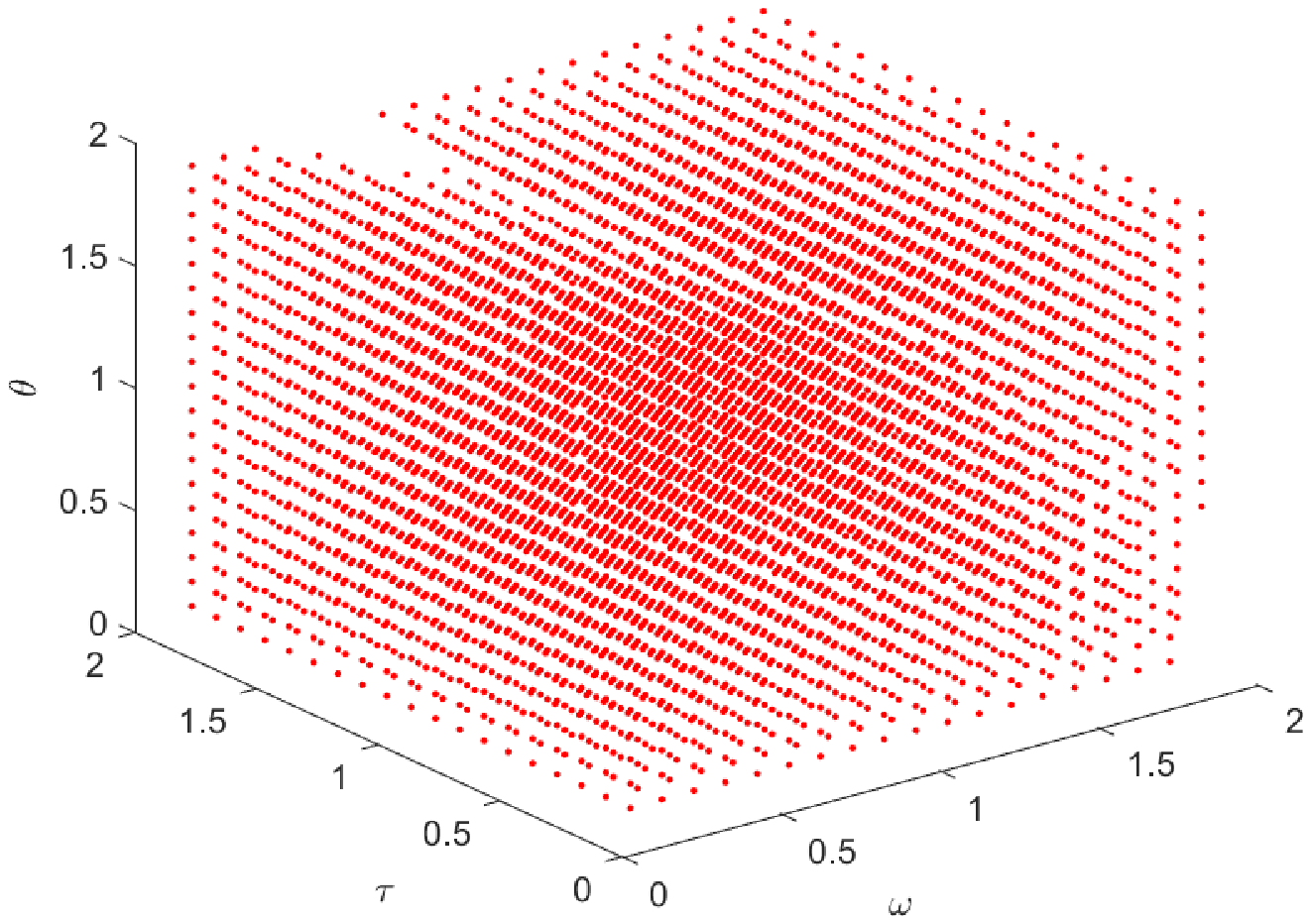}
    \hfill
    \includegraphics[width=0.4\textwidth]{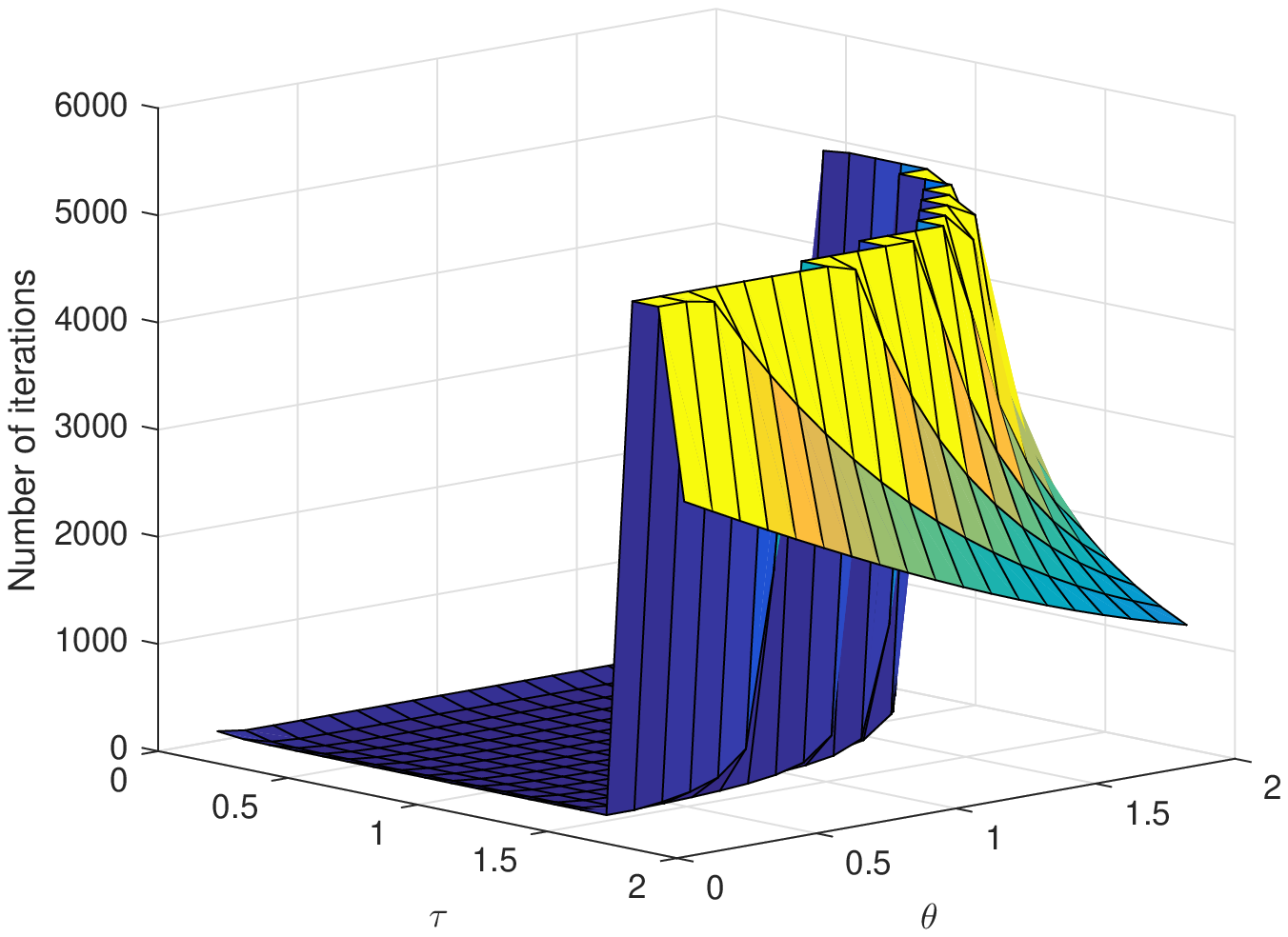}
    \\
    \includegraphics[width=0.4\textwidth]{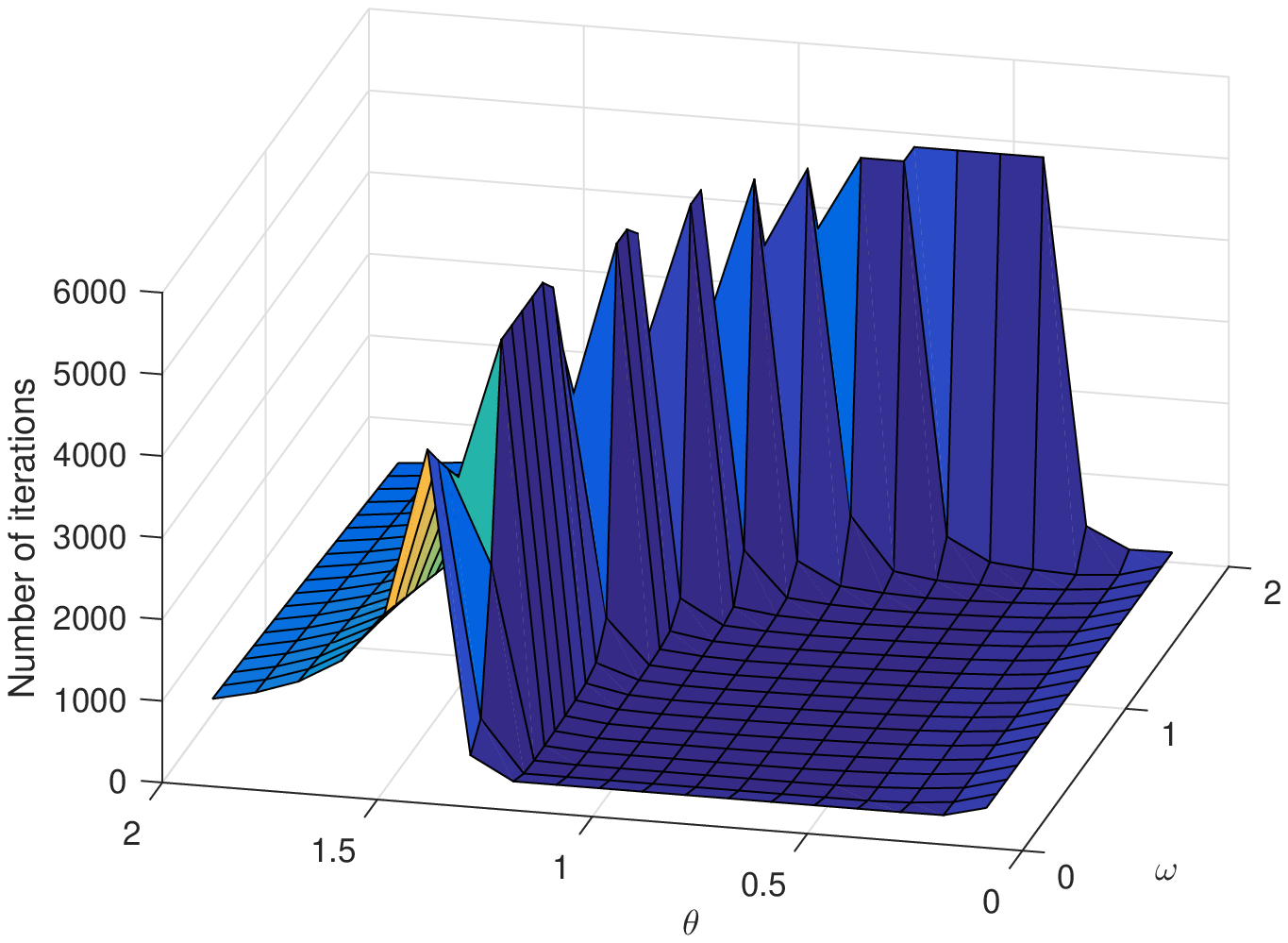}
    \hfill
    \includegraphics[width=0.4\textwidth]{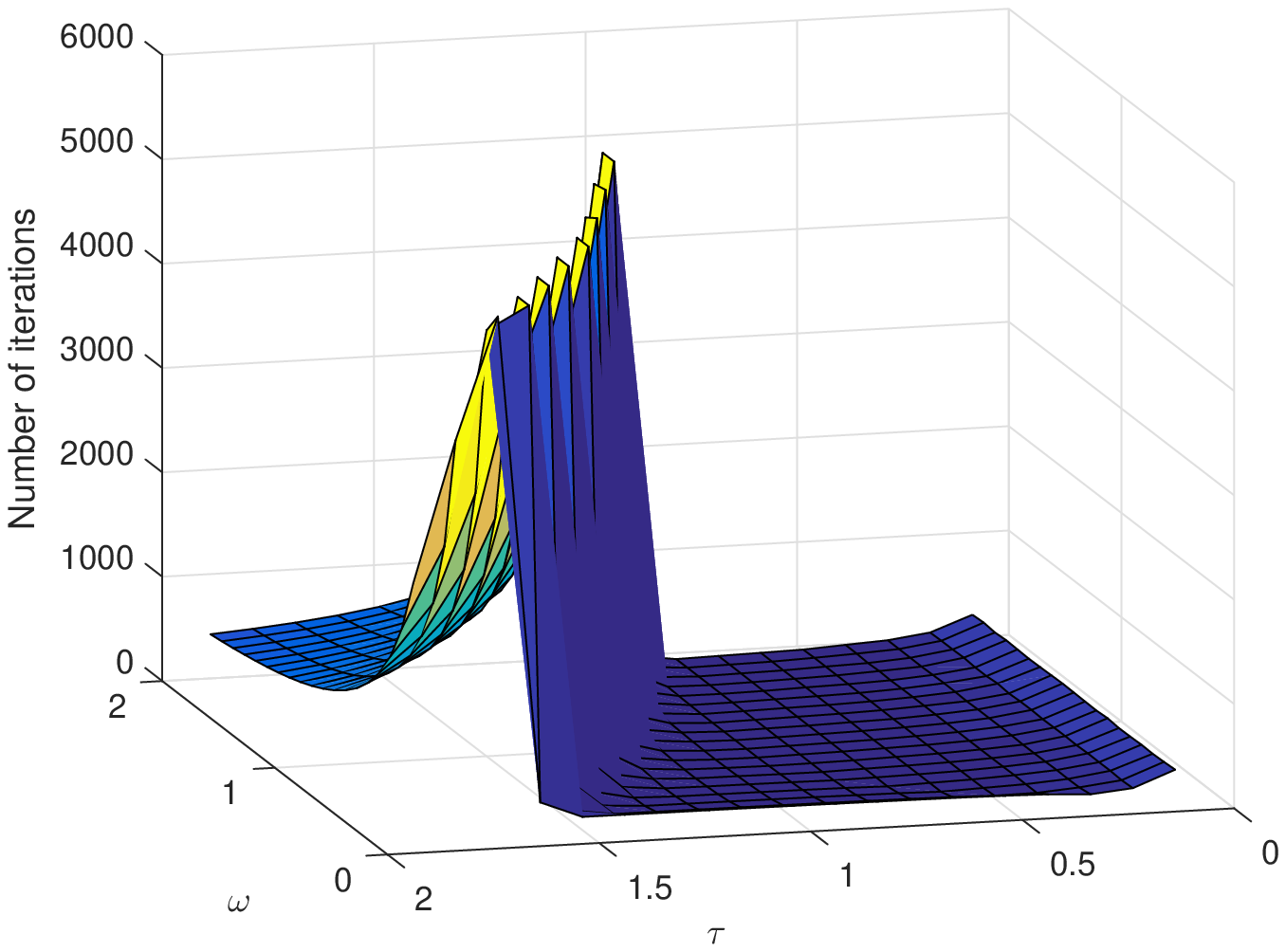}
    \caption{\label{fig:plots}%
        Top left: The region of parameter values for which GSOR satisfies ${\rm Res}\le 10^{-8}$ within $5,000$ iterations.
	    Other plots: Characteristic curves for the number of iterations versus parameters $\omega$, $\tau$ and $\theta$ for GSOR with $\omega=1$ (top right), $\tau = 1$ (bottom left), and $\theta = 1$ (bottom right).
	    All plots are for saddle-point systems from the liquid crystal directors model with $n=3069$, $m=p=1023$.
    }
\end{figure}

\begin{figure}[ht]   
	\centering
	\subfloat[$\calA$]{
		\includegraphics[width=0.4\linewidth]{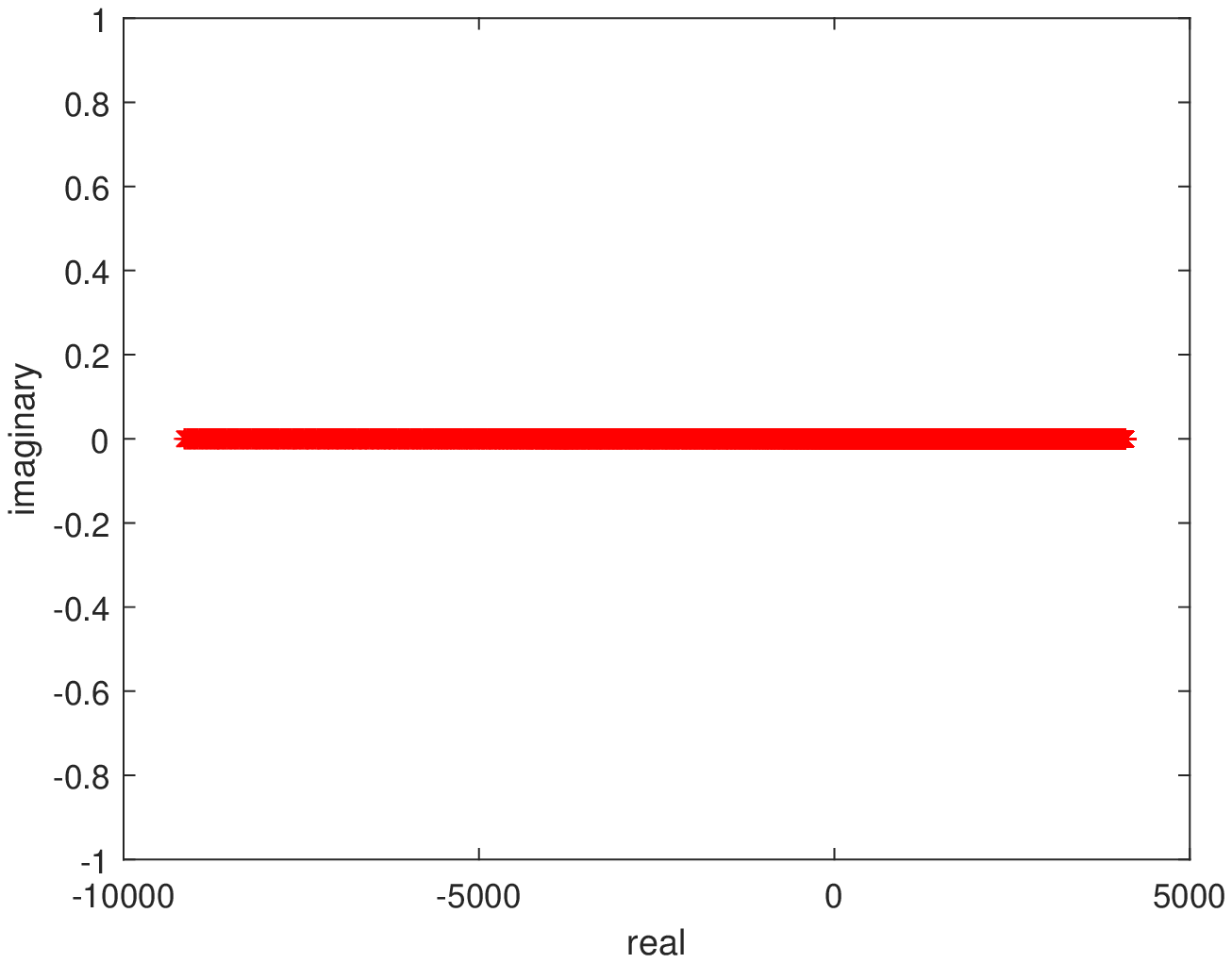}
	}
	\subfloat[${\calP}^{-1}\calA$ with $\tau=0.1,\theta=1$]{
		\includegraphics[width=0.4\linewidth]{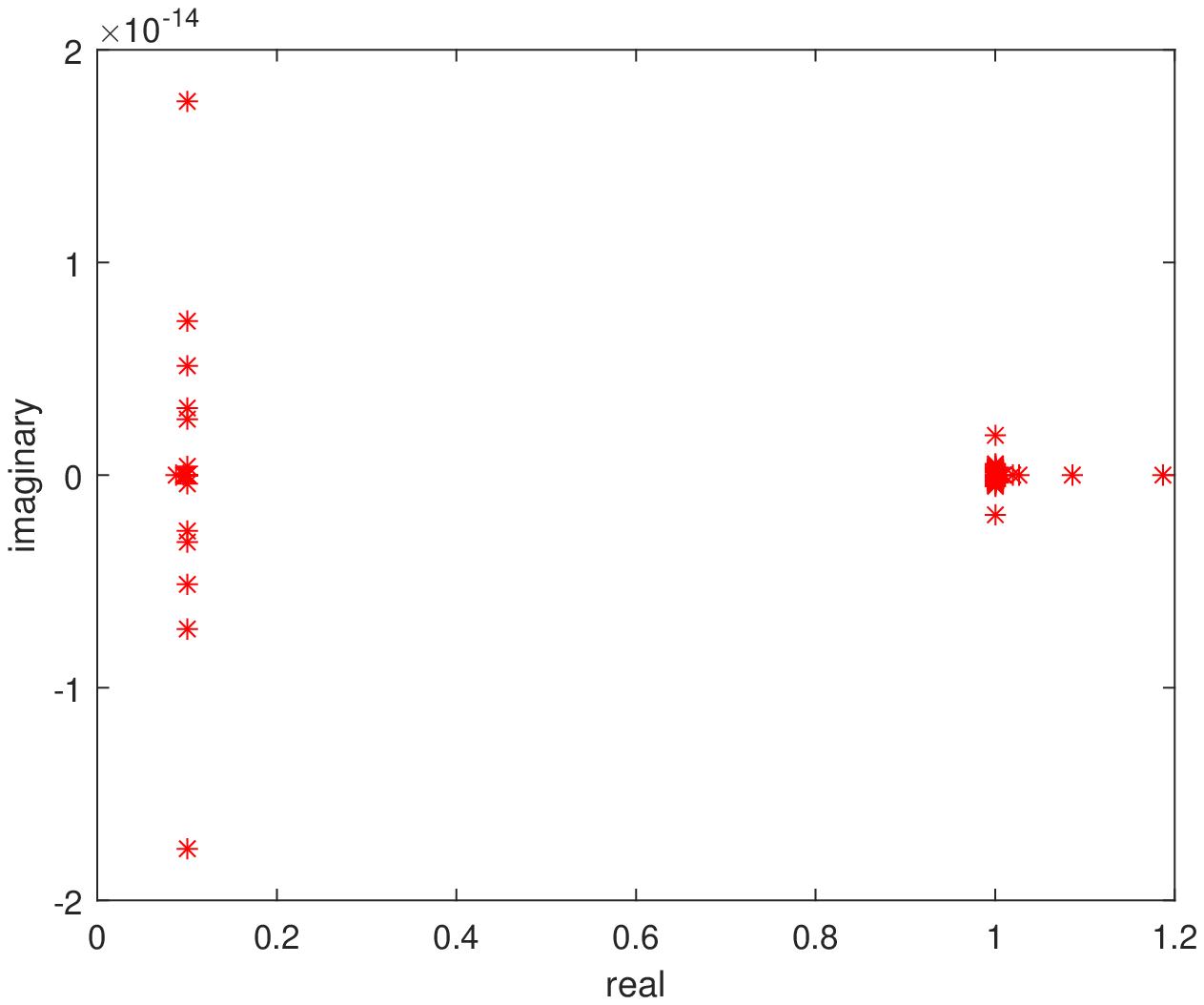}
	}\\
	\subfloat[${\calP}^{-1}\calA$ with $\tau=1,\theta=1$]{
		\includegraphics[width=0.4\linewidth]{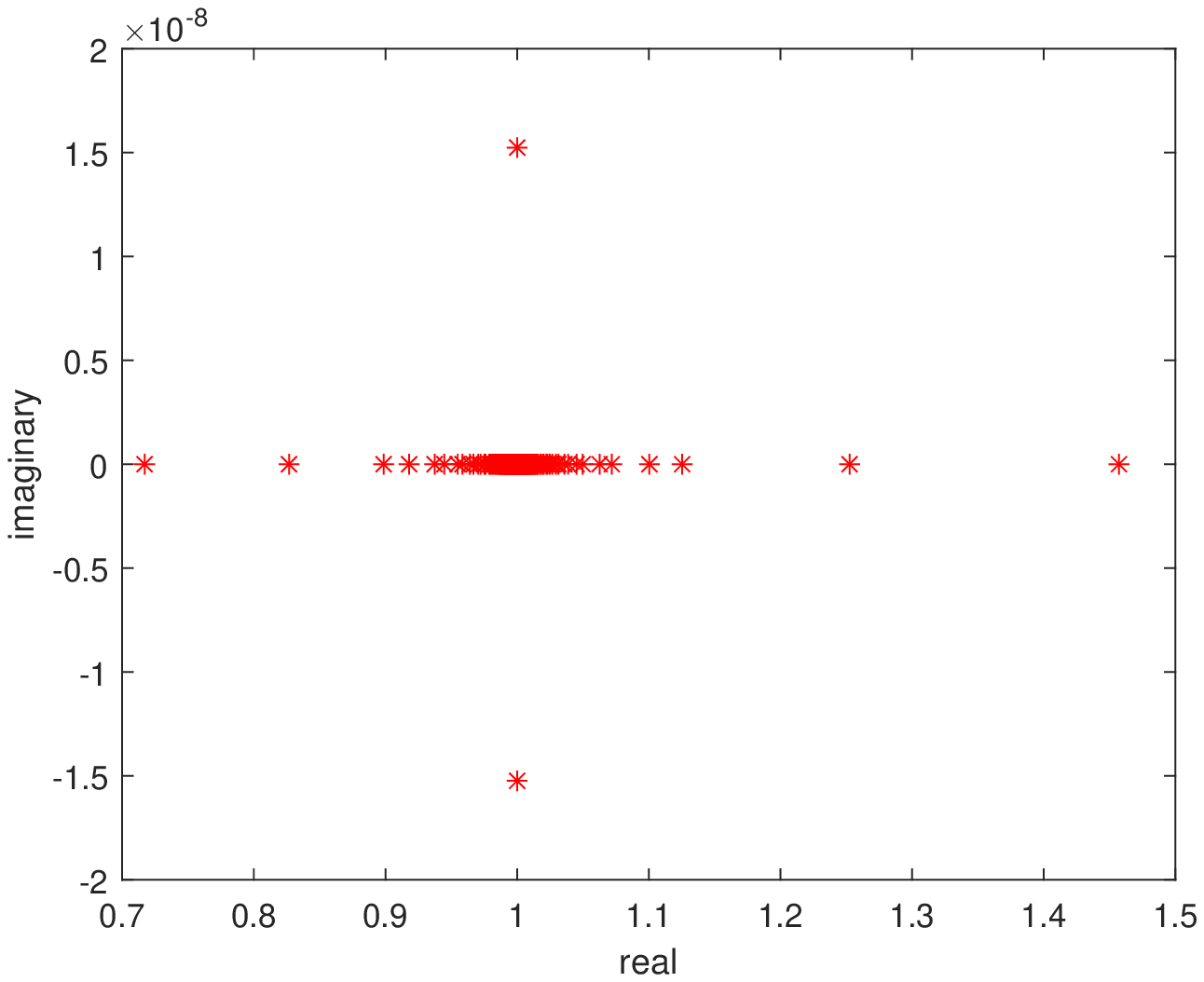}
	}
	\subfloat[${\calP}^{-1}\calA$ with $\tau=1,\theta=0.1$]{
		\includegraphics[width=0.4\linewidth]{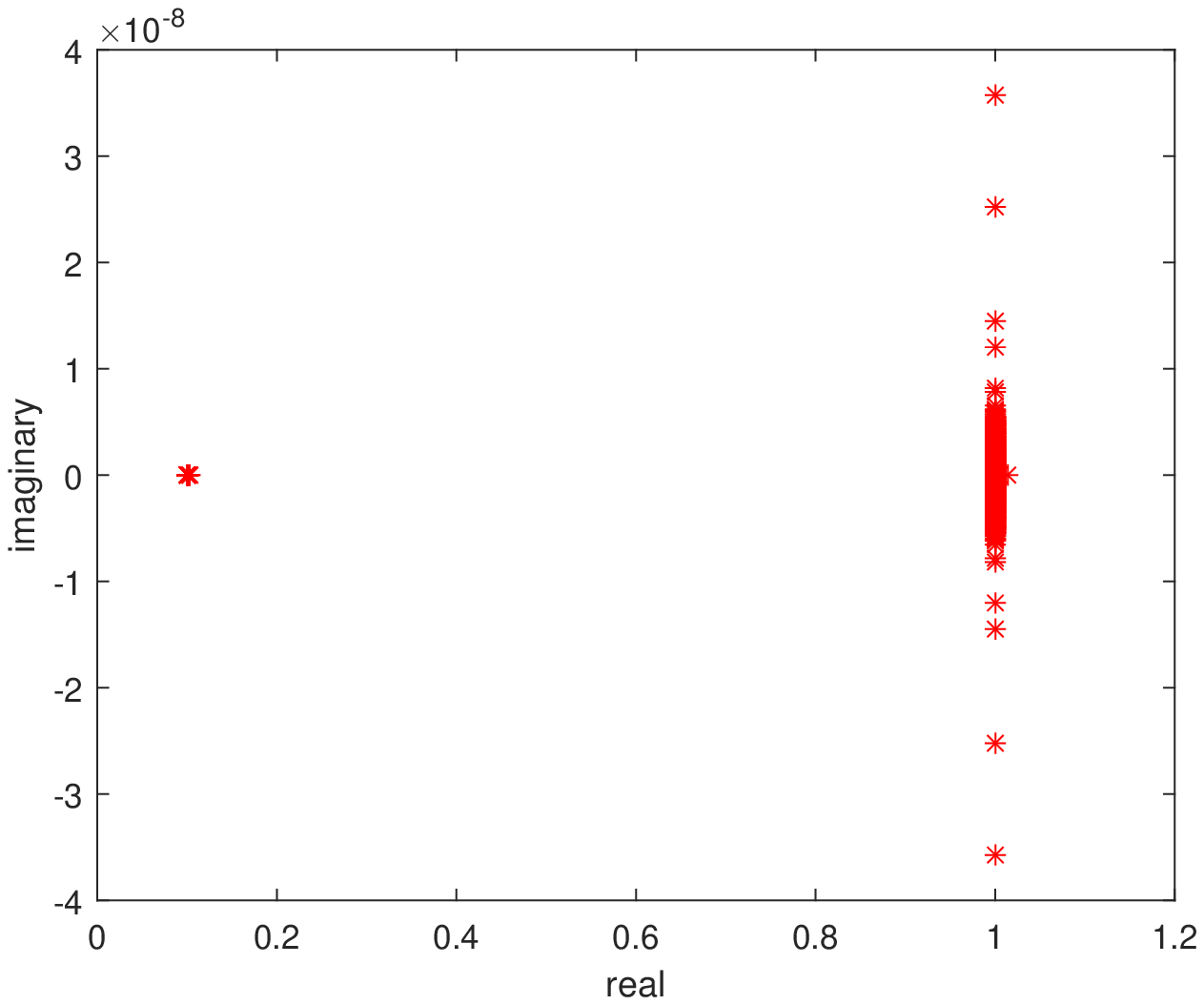}
	}
	\caption{Eigenvalue distributions of the original matrix and the GSOR preconditioned matrices for saddle-point systems from the liquid crystal directors model with $n=3069,m=p=1023$.}
	\label{fig:eigenvaluelsdm}
\end{figure}

To see the role of the parameters in the convergence behavior of GSOR, \cref{fig:plots} shows the region of the parameters where GSOR satisfies ${\rm Res}\le 10^{-8}$ within $5,000$ iterations, and the characteristic curves of the number of iterations versus the parameters for $N=1,023$. In \cref{fig:eigenvaluelsdm}, we plot the eigenvalue distributions of the original matrix and the GSOR preconditioned matrix ${\calP}^{-1}{\calA}$ with different $\tau$ and $\omega$.


\subsection{Saddle-point systems from the mixed Stokes-Darcy model in porous media applications}\label{subsec:msdm}

Fluid flow in $\Omega_{f}\subset \R^2$ coupled with porous media flow in $\Omega_{p}\subset \R^2$ is governed by the static Stokes equations
\begin{equation}\label{sdm}
      -\upsilon\Delta\,{\bm u}_f + \nabla\, p_f  = {\bm f},
      \quad \textup{and} \quad
      {\rm div}\,{\bm u}_f = 0,
      \quad {\bm x}\in \Omega_{f},
\end{equation}
where $\Omega_{f}\cap \Omega_{p}=\emptyset$ and $\overline{\Omega}_{f}\cap \overline{\Omega}_{p}=\Gamma$ with $\Gamma$ being an interface, $\upsilon>0$ is the kinematic viscosity, and $\bm f$ is the external force.


In the porous media region, the governing variable is $\phi = \frac{p_p}{\rho_f g}$, where $p_p$ is the pressure in $\Omega_{p}$, $\rho_f$ is the fluid density, and $g$ is the gravity acceleration. The velocity ${\bm u}_p$ of the porous media flow is related to $\phi$ by Darcy's law and is also divergence free:
\begin{equation}\label{sdm2}
      {\bm u}_p = -\dfrac{\epsilon^2}{r\upsilon}\nabla \phi
      \quad \textup{and} \quad
      -{\rm div}\,{\bm u}_p = 0,
      \quad {\bm x}\in \Omega_{p},
\end{equation}
where $r$ is the volumetric porosity, and $\epsilon$ the characteristic length of the porous media.

Applying finite element discretization to the mixed Stokes-Darcy model~\eqref{sdm}--\eqref{sdm2} with the Dirichlet boundary conditions leads to linear systems of form~\eqref{a1}~\cite{Cai2009}.

In our numerical experiments, we set $\upsilon=1$, $r=1$, and $\epsilon=\sqrt{0.1}$. The computational domain is $\Omega_f=(0,1)\times(1,2)$, $\Omega_p=(0,1)\times(0,1)$ and the interface is $\Gamma=(0,1)\times\{1\}$. We use a uniform mesh with grid parameters $h=2^{-3},\,2^{-4},\,2^{-5},\,2^{-6}$ to decompose $\Omega_f$, P2--P1 elements in the fluid region, and P2 Lagrange elements in the porous media region.

For this problem, $P$ is the pressure mass matrix discretized from the decoupled problem of~\eqref{sdm}--\eqref{sdm2}~\cite{Cai2009}.  In all algorithms we use Cholesky factorization to solve the systems $Ax=r$, $Py=r$ and $BA^{-1}B\T y=r$.
Numerical results for saddle-point systems from the mixed Stokes-Darcy model~\eqref{sdm}--\eqref{sdm2} are listed in \Cref{tab21,tab22}, where the parameters choices for GSOR and the corresponding notation are as follows.

\begin{center}\scriptsize
   \medskip
   \begin{tabular}{|c |c |c |c |c|}\hline
     Method & GSORa & GSORb & GSORc & GSORd
   \\ \hline
    $(\omega,\tau,\theta)$  & $(0.5,1.5,1.0)$ & $(0.5, 1.7, 0.8)$& $(0.5,1.6,1.2)$ &$(0.6,1.5,1.0)$
   \\ \hline
   \end{tabular}
   \medskip
\end{center}

For this problem, $\nu_{\max}=1.0057$. The matrix $A-C\T D^{-1}C$ is no longer SPD, so convergence of the Uzawa-like method cannot be guaranteed \cite[Theorem 3]{Benzi2019}. We tested several $\alpha$ ranging from $0.005$ to $0.5$ for $h=2^{-3}$. Uzawa failed in all cases. Thus, we do not report results for Uzawa in \cref{tab22}. As MINRES and GMRES worked only for systems with $h \ge 2^{-5}$, we again do not report their results.

\begin{figure}[ht]   
	\centering
	\includegraphics[width=0.4\linewidth]{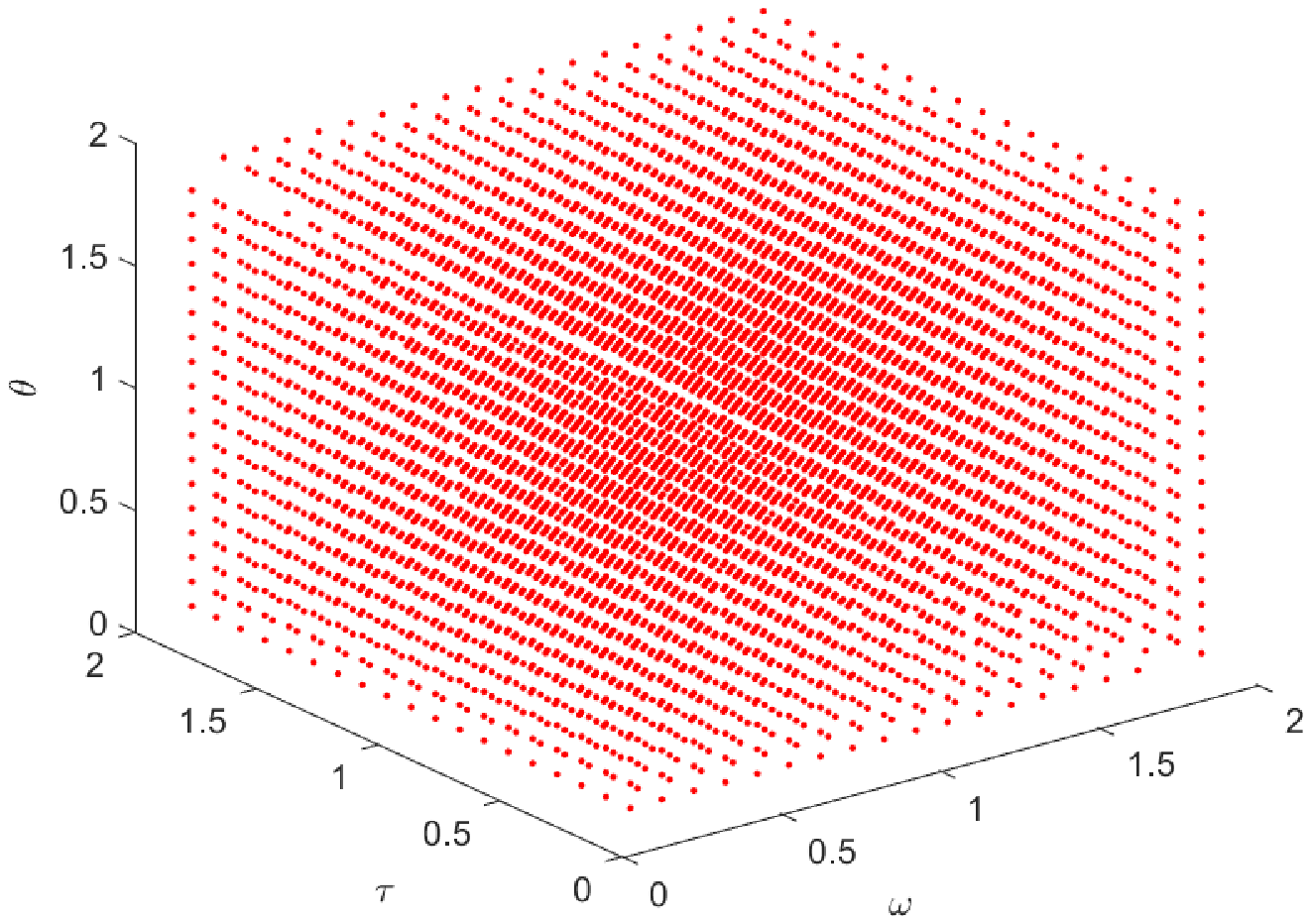}
	\hfill
	\includegraphics[width=0.4\linewidth]{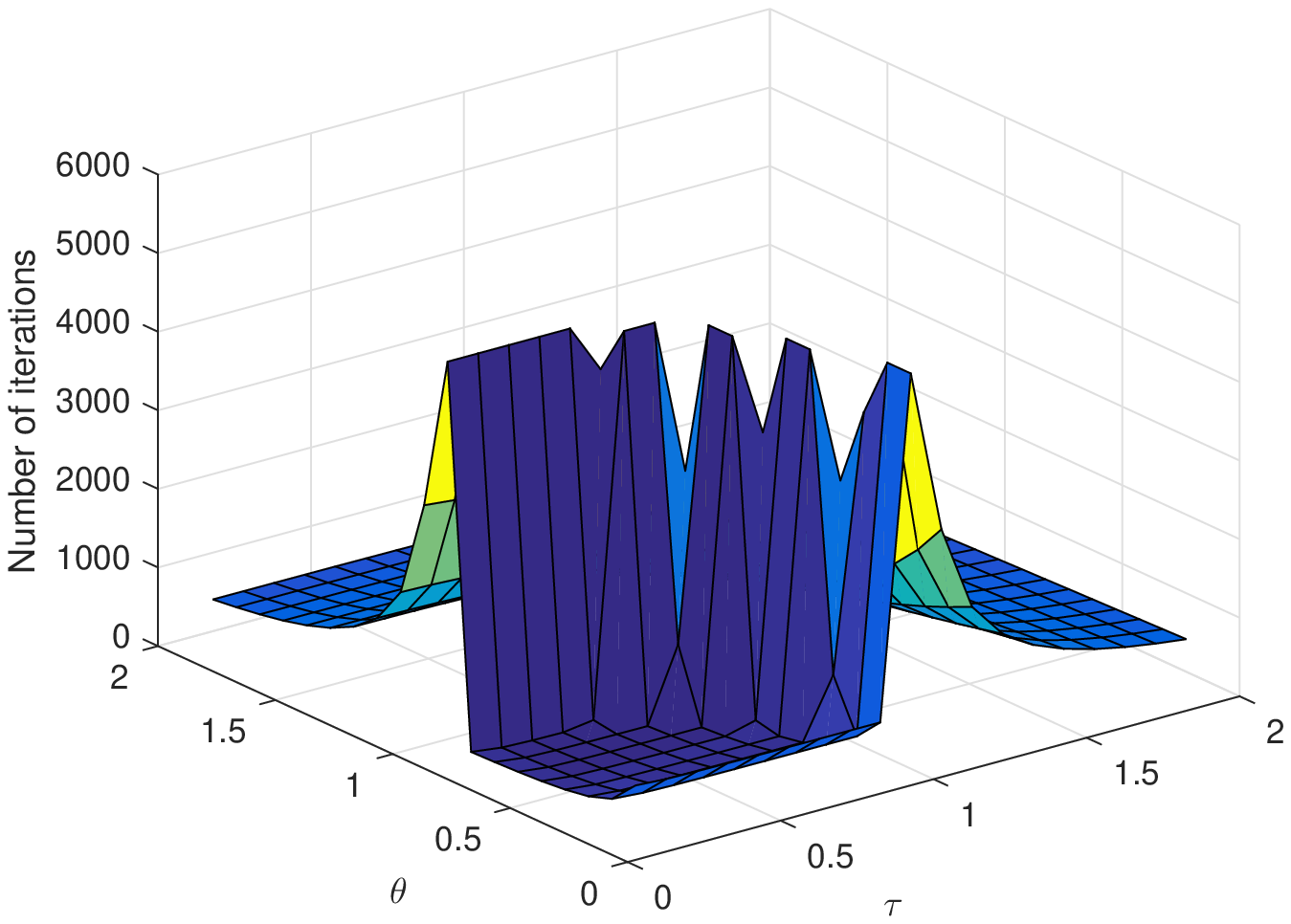}
	\\
	\includegraphics[width=0.4\linewidth]{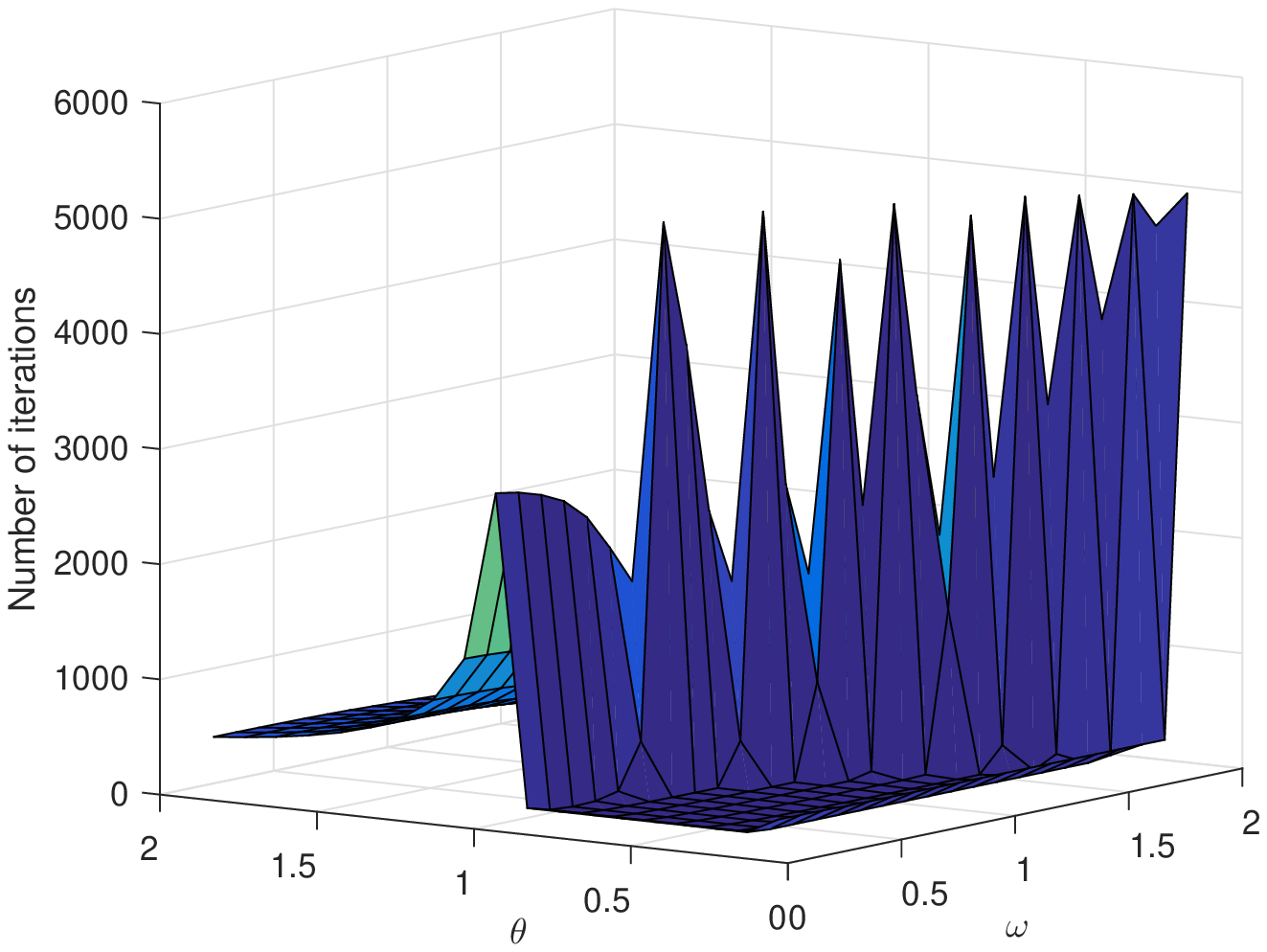}
	\hfill
	\includegraphics[width=0.4\linewidth]{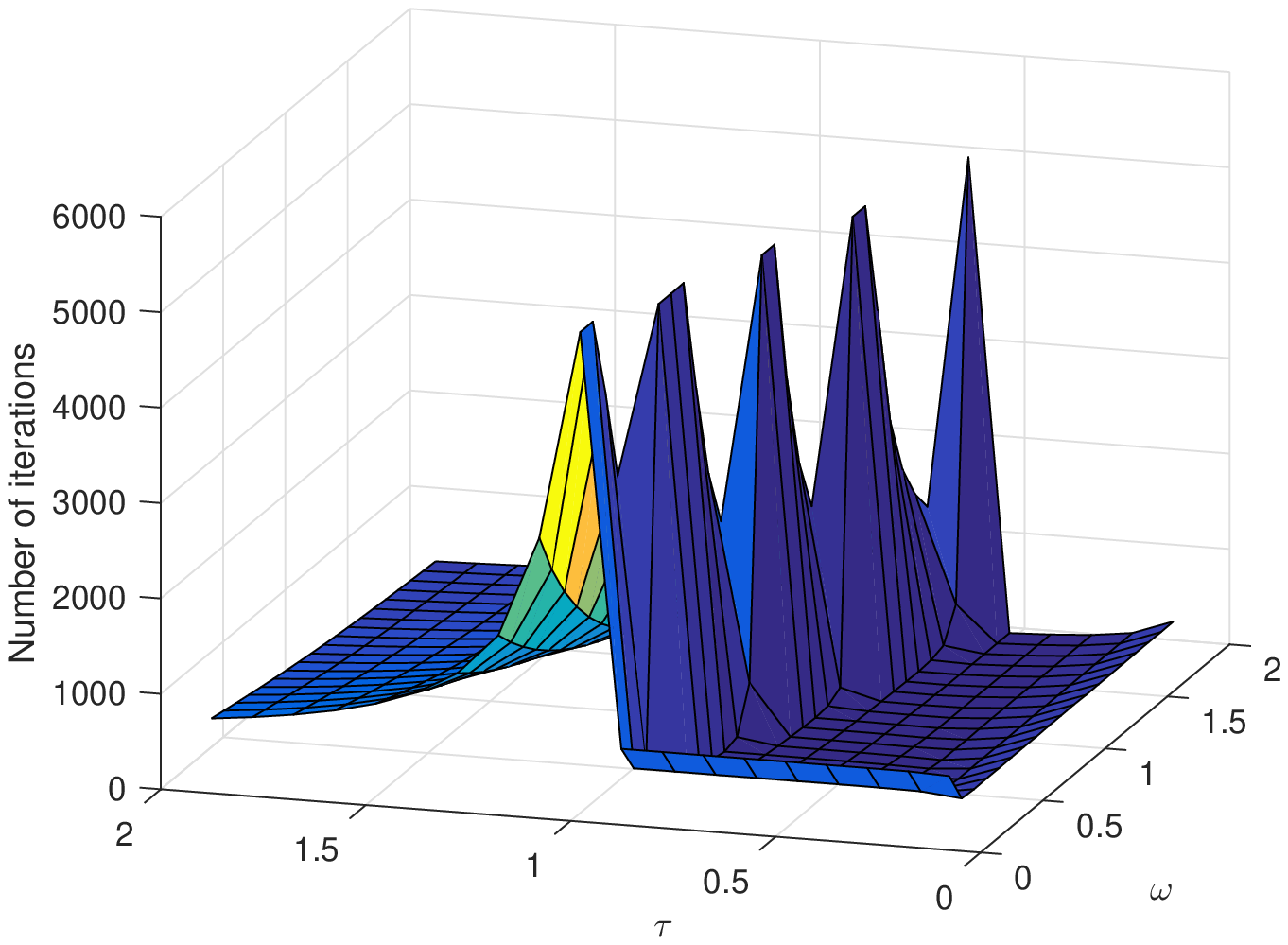}
	\caption{\label{fig:st-plots}%
	    Top left: The region of parameter values for which GSOR satisfies ${\rm Res}\le 10^{-8}$ within $5,000$ iterations.
	    Other plots: Characteristic curves for the number of iterations versus parameters $\omega$, $\tau$ and $\theta$ for GSOR with $\omega=1$ (top right), $\tau = 1$ (bottom left), and $\theta = 1$ (bottom right).
	    All plots are for saddle-point systems from the mixed Stokes-Darcy model with $n=578$, $m=81$, $p=289$.
	}
\end{figure}

\begin{figure}[ht]   
	\centering	
	\subfloat[$\calA$]{%
		\includegraphics[width=.4\linewidth]{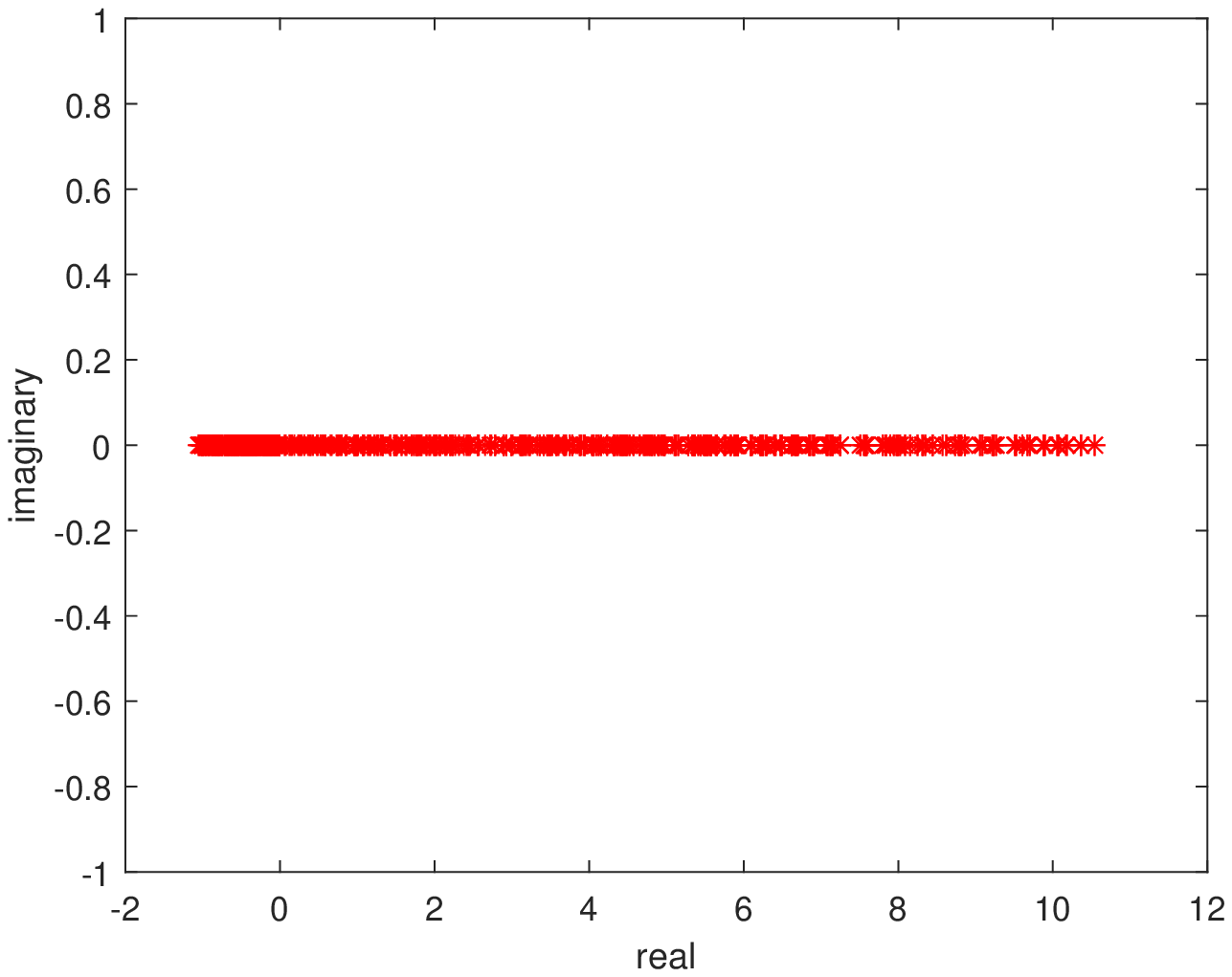}%
	}
	\subfloat[${\calP}^{-1}\calA$ with $\tau=0.1,\theta=1$]{%
		\includegraphics[width=.4\linewidth]{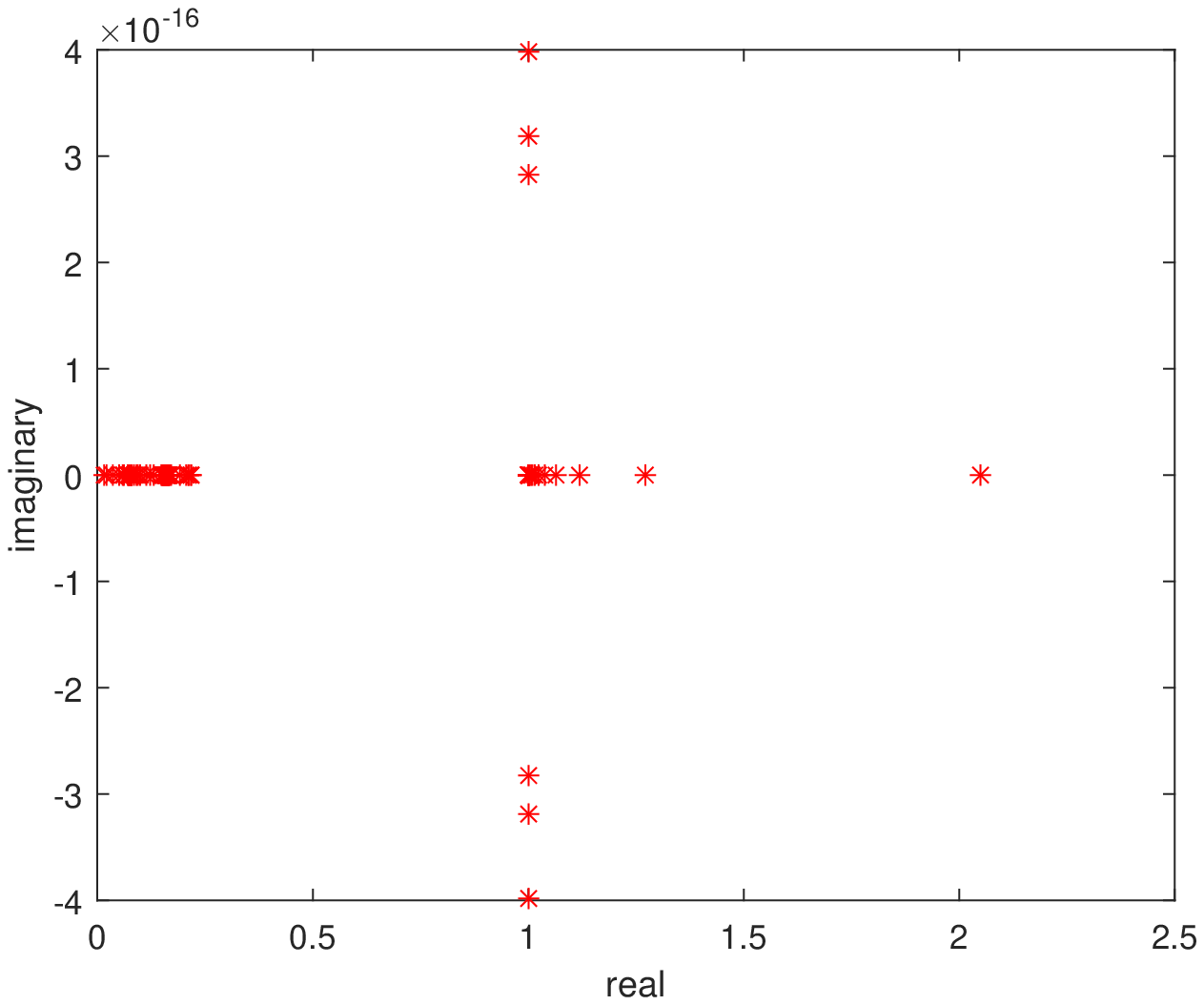}%
	}\\
	\subfloat[${\calP}^{-1}\calA$ with $\tau=1,\theta=1$]{%
		\includegraphics[width=.4\linewidth]{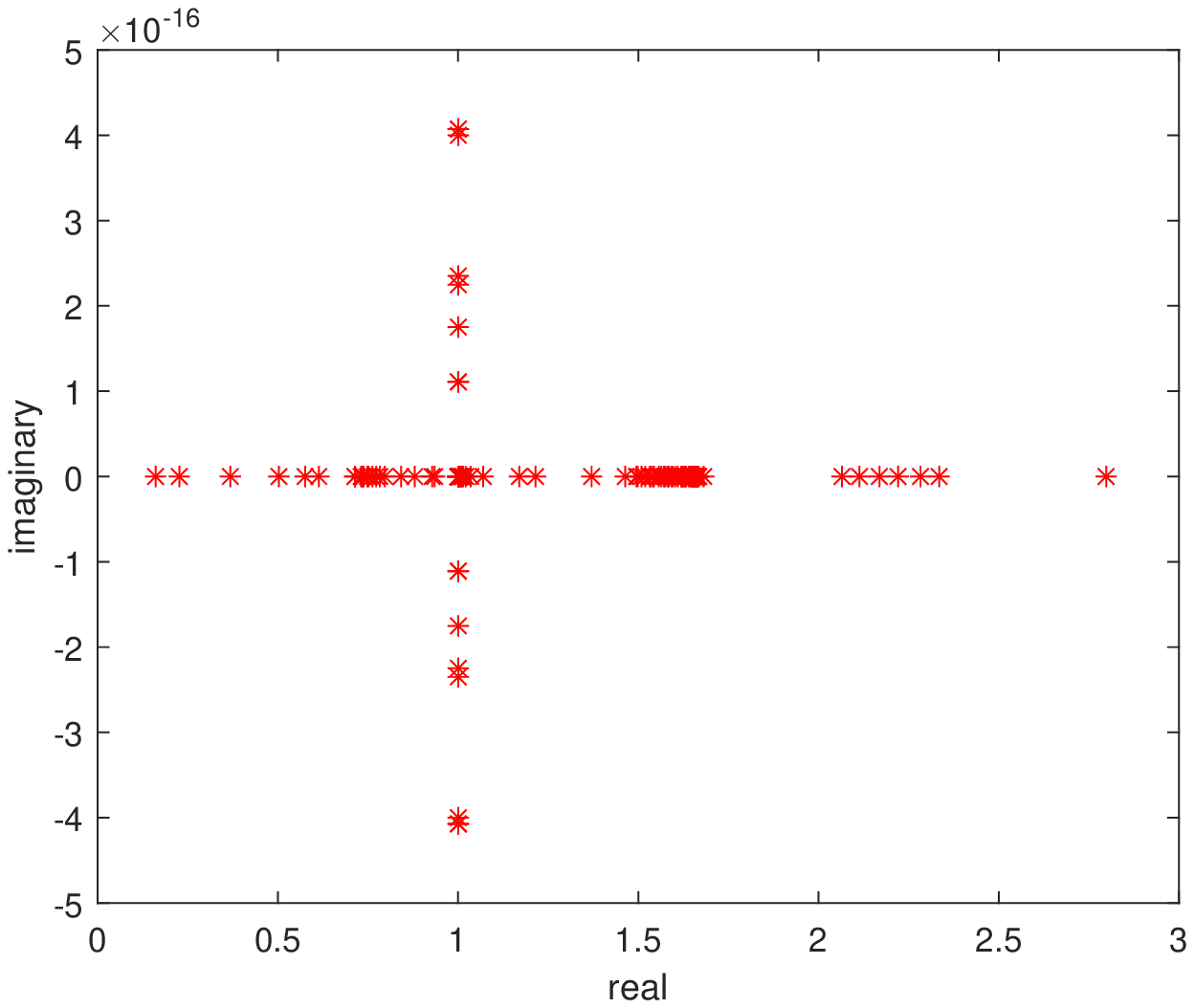}%
	}
	\subfloat[${\calP}^{-1}\calA$ with $\tau=1,\theta=0.1$]{%
		\includegraphics[width=.4\linewidth]{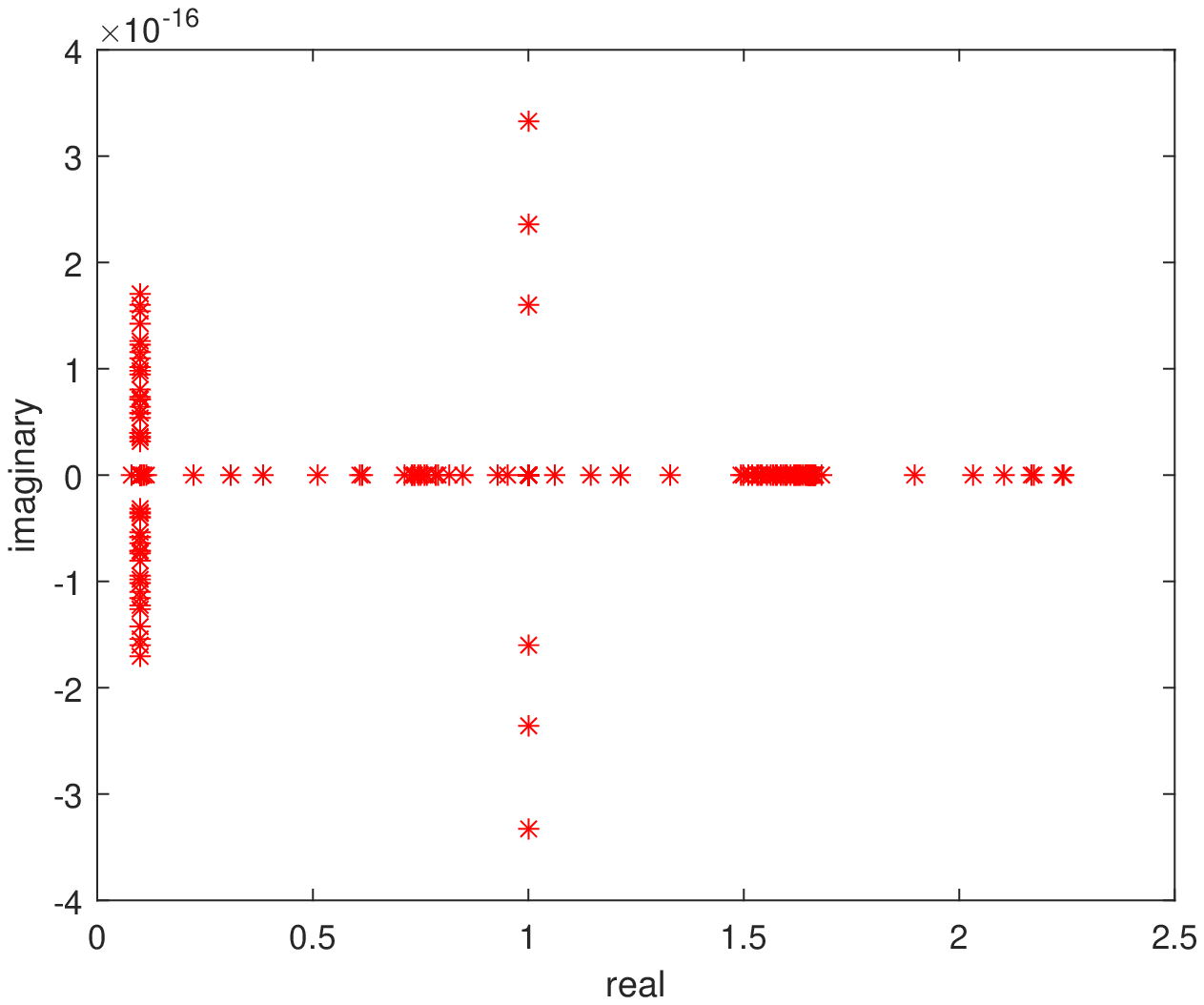}%
	}
	\caption{Eigenvalue distributions of the original matrix and the GSOR preconditioned matrices for saddle-point systems from the mixed Stokes-Darcy model with $n=578,m=81,p=289$.}
	\label{fig:eigenvaluesd}
\end{figure}

To see the role of the parameters in the convergence behavior of GSOR, \Cref{fig:st-plots} shows the region of parameters for which GSOR satisfies ${\rm Res}\le 10^{-8}$ within $5,000$ steps, and the characteristic curves of iteration numbers versus parameters for $h=2^{-3}$. In \cref{fig:eigenvaluesd}, we plot the eigenvalue distributions of the original matrix and the GSOR preconditioned matrix ${\calP}^{-1}{\calA}$ with different $\tau$ and $\omega$.


\begin{table}[ht]\scriptsize   
  \centering
  \caption{The CPU time of the Cholesky factorization.}\label{tab21}
  \begin{tabular}{|c |c |c |c |c |c |c |c |c |}
    \hline
$h$& $2^{-3}$&$2^{-4}$&$2^{-5}$&$2^{-6}$&$2^{-7}$\\\hline
$n$ & 578 & 2178 & 8450 & 33282 & 132098 \\\hline
$m$ & 81 & 289 & 1089 & 4225 & 16641\\\hline
$p$ & 289 & 1089 & 4225 & 16641 & 66049 \\\hline
$n+m+p$ & 948 & 3556 & 13764 &54148& 214788 \\\hline
$A$ &  0.0008& 0.0048&0.029  & 0.23&1.75  \\\hline
$P$ & 0.0003 &0.0004  &0.0011  & 0.02& 0.10\\\hline
$BA^{-1}B\T$ &0.0064 & 0.18   & 7.18  & 555.59 &31368.15 \\\hline

  \end{tabular}
\end{table}

\begin{table}[ht]\scriptsize   
  \centering
  \caption{Numerical results for saddle-point systems from mixed the Stokes-Darcy model.}\label{tab22}
  \begin{tabular}{|c|c|c|c|c|c|c|c|c|c|c|}
    \hline
    & $h$& $2^{-3}$&$2^{-4}$&$2^{-5}$&$2^{-6}$&$2^{-7}$ \\\hline
      &Iter  & 50    &   49 & 49 &47 & 47  \\
GSORa &CPU   &0.05 & 0.15 &0.89   & 10.19& 54.22\\
      &Res   &5.99e-09&9.80e-09 & 8.41e-09 & 9.53e-09 &6.85e-09\\\hline

      &Iter  &50   &    50&  50  & 50 & 50\\
GSORb &CPU   &0.03 &  0.15 & 0.91  &  10.13& 64.10 \\
      &Res   &6.54e-09 &5.93e-09 &5.51e-09   &5.82e-09 &5.74e-09\\\hline

      &Iter  & 50      & 50 & 50   & 50 & 49 \\
GSORc &CPU   & 0.04 & 0.15   &0.90   &10.02  & 62.91 \\
      &Res   & 5.98e-09 & 5.19e-09& 6.92e-09  & 7.00e-09&9.53e-09 \\\hline

      &Iter  & 48     &   45 & 42   &  39 & 38\\
GSORd &CPU   & 0.04   &  0.14&  0.80   & 7.70  &49.98 \\
      &Res   &8.39e-09& 8.53e-09 & 8.52e-09     &9.63e-09 &5.26e-09 \\\hline

      &Iter  &  158    &  150&  141    & 132 & 124\\
GBSORa &CPU   &0.21  &0.84   &   6.01  &  90.33  & 948.46  \\
      &Res   & 9.79e-09& 9.16e-09&  9.46e-09   &  9.96e-09 &9.77e-09\\\hline

      &Iter  & 73   &69 &    65   & 61 & 58\\
GBSORb &CPU   & 0.08   & 0.36  &  2.74  &  41.47 & 441.51  \\
      &Res   &  9.17e-09 &9.18e-09&  9.28e-09   &  9.57e-09& 8.28e-09\\\hline

      &Iter  & 44   &42  &   40  &  37 & 35\\
GBSORc &CPU   & 0.04 & 0.19  &  1.73    & 25.59   & 265.67 \\
      &Res   &  9.33e-09 & 8.25e-09  &  7.38e-09 & 9.57e-09& 8.61e-09 \\\hline

      &Iter  & 767.5  & 1491 &  2997.5   & 5912.5 & 13826.5\\
BICGSTAB &CPU   &  0.09  &0.41   &  2.23  & 21.51 & 288.07 \\
      &Res   & 7.96e-09 &7.64e-09 &  9.57e-09   & 4.76e-09 & 9.48e-09 \\\hline

      &Iter  & 18 &  18 &  18   & 17  &  18\\
BPMINRES &CPU   & 0.13 & 0.71 &  6.11   & 216.39 & 3479.83\\
      &Res   &8.34e-09& 3.66e-09&   1.49e-09   & 5.70e-09 &2.85e-09\\\hline

      &Iter  & 10  &  10    &  10   &  10  &10 \\
BPGMRES &CPU   & 0.14 & 0.81  &   10.00  & 323.95  & 3992.25 \\
      &Res   & 1.60e-09 &  1.56e-09&  1.05e-09  &  7.06e-10  &1.39e-09 \\\hline

      &Iter  & 24  &  24    &25  &  26  &27  \\
GPGMRES &CPU   &0.14  &0.84 &   3.53  &17.34   & 92.48\\
      &Res   & 1.06e-09 &  4.98e-09 & 6.87e-09  & 1.53e-09  &6.15e-10 \\\hline
  \end{tabular}
\end{table}

\Cref{tab11,tab12,tab21,tab22} and \Cref{fig:plots,fig:eigenvaluelsdm,fig:st-plots,fig:eigenvaluesd} illustrate that GSOR is a practical method,  and its advantages increase with the problem size. We see from \Cref{tab11,tab12,tab21,tab22} that BPMINRES and BPGMRES are not practical in terms of CPU times.  \Cref{fig:plots,fig:st-plots} indicate that the convergence rate of GSOR depends strongly on $\omega$, $\tau$ and $\theta$. \Cref{fig:eigenvaluelsdm,fig:eigenvaluesd} show that $\calP$ greatly improves the eigenvalue distribution of the original $\calA$.


\section{Conclusions}\label{sec:con}

We presented a theoretical and numerical study of the GSOR method for solving the double saddle-point problem~\eqref{a1}. GSOR is convergent with suitable parameters $\omega$, $\tau$, and $\theta$. Unlike existing work, our proof is based on the necessary and sufficient conditions for all roots of a real cubic polynomial to have modulus less than one.
We analyzed a class of block lower triangular preconditioners $\calP$
induced from GSOR and derived explicit and sharp bounds for the
eigenvalues of preconditioned matrices.
The numerical results presented are highly encouraging. GSOR requires the least CPU time, and especially for larger problems, its advantages are clear. A shortcoming is the need to choose the three parameters. 
A practical method to choose them is a topic for future research.

\clearpage

\section*{Acknowledgments}

This research began with the work of our colleague and friend Dr Oleg Burdakov. 
In 2019, Oleg focused on Barzilai-Borwein-type methods to solve quasi-definite linear systems, and he conducted preliminary tests on double saddle-point problems. While testing the BB-type methods, we found that GSOR performs well on the double saddle-point problems and were motivated to start this work.  We are grateful for Oleg's insight and foresight.

We are also grateful to Mingchao Cai, Alison Ramage, and Zhaozheng Liang for providing the test problems used in our numerical experiments.

\bibliographystyle{siamplain.bst}
\bibliography{references}

\end{document}